\documentclass[12pt,a4paper]{amsart}
\usepackage[english]{babel}
\usepackage[T1]{fontenc}
\usepackage{mathtools}
\usepackage{fancyhdr}
\usepackage{appendix}
\usepackage{amssymb,amscd,amsxtra,calc}
\usepackage{mathrsfs}
\usepackage{xcolor}
\usepackage{multirow}
\usepackage[all]{xy}
\usepackage{longtable}
\usepackage[colorlinks=true,linkcolor=blue,anchorcolor=blue,citecolor=blue,pagebackref,linktocpage]{hyperref}
\usepackage{cleveref}
\usepackage{tikz}
\usepackage{cite}
\newcommand{\RzeroTag}{$\spadesuit$}
\newtheorem{theorem}{Theorem}[section]
\newtheorem{proposition}[theorem]{Proposition}
\newtheorem{lemma}[theorem]{Lemma}
\newtheorem{corollary}[theorem]{Corollary}
\newtheorem{example}[theorem]{Example}
\newtheorem{definition}[theorem]{Definition}
\newtheorem{remark}[theorem]{Remark}
\newtheorem{question}[theorem]{Question}
\newtheorem{problem}[theorem]{Problem}
\newcommand{\C}{\mathbb C}
\newcommand{\Z}{\mathbb Z}
\newcommand{\R}{\mathbb R}
\newcommand{\PP}{\mathbb P}
\newcommand{\OO}{\mathcal O}
\newcommand{\Pic}{\operatorname{Pic}}

\newcommand{\MW}{\operatorname{MW}}
\newcommand{\Cox}{\operatorname{Cox}}
\newcommand{\Eff}{\operatorname{Eff}}
\newcommand{\Nef}{\operatorname{Nef}}
\newcommand{\NE}{\overline{\operatorname{NE}}}
\newcommand{\gr}{\operatorname{gr}}
\newcommand{\ord}{\operatorname{ord}}

\newcommand{\arxiv}[1]{\href{https://arxiv.org/abs/#1}{arXiv:#1}}
\newcommand{\pos}[1]{\max\{#1,0\}}
\allowdisplaybreaks

\title[Cox rings of Jacobian elliptic surfaces]{Cox rings of Jacobian elliptic surfaces of Kodaira dimension one}
\author{Antonio Laface}
\address{Departamento de Matem\'atica, Universidad de Concepci\'on, Casilla 160-C, Concepci\'on, Chile}
\email{\href{mailto:alaface@udec.cl}{alaface@udec.cl}}

\author{Sichen Li}
\address{
School of Mathematics, East China University of Science and Technology, Shanghai 200237, P. R. China}
\email{\href{mailto:sichenli@ecust.edu.cn}{sichenli@ecust.edu.cn}}

\begin{document}
\begin{abstract}
We study Cox rings of relatively minimal semistable Jacobian elliptic
surfaces over $\mathbb P^1$ with $\chi=\chi(\mathcal O_X)\ge3$ and finite
Mordell--Weil group. When the Mori cone is generated by the zero section
and the reducible-fiber components, we give a finite criterion for Cox-ring
generation: a numerically determined set of nef classes must admit sections
nonvanishing along the zero section. The criterion supplies a minimal
homogeneous generating set without assuming finite generation in advance.
We prove that its nonvanishing condition holds whenever
$\sum_i n_i\le2\chi+4$, where $I_{n_i}$ are the reducible fibers; hence every
such surface is a Mori dream surface. Independently of this condition,
we exclude ample generator degrees throughout the stated cone range.
We also determine the section--fiber subring and give minimal
complete-intersection presentations when all $s$ reducible fibers have
type $I_2$ and $s\le2\chi+1$, including the endpoint beyond that cone range.
\end{abstract}
\maketitle
\tableofcontents
\section{Introduction}
The classification of surfaces with finitely generated Cox ring remains
open. Explicit computations and geometric criteria are available for
rational surfaces, K3 surfaces, and several classes of surfaces of general
type; see \cite{ADL21,ADR23,AG25,AGL16,AHL10,AL11,KL19,KL26,TVV11}.
For Jacobian elliptic surfaces of Kodaira dimension one, the cone and
semiampleness results of \cite{LLL26,LL26} provide a starting point.
A rational polyhedral Mori cone alone, however, does not ensure finite
generation of the Cox ring. The remaining geometry is reflected in the
sections of nef divisors on its boundary.

Our main result makes this issue finite: under the cone hypothesis below,
we construct a finite set of nef classes and show that nonvanishing along
the zero section in these classes produces a minimal homogeneous
generating set of the entire Cox ring. We then establish this condition
in the numerical range $\sum_i n_i\le2\chi+4$. Alongside this criterion,
we give explicit presentations for the section--fiber subring and for
configurations of $I_2$ fibers, and a vanishing argument excluding all
ample generator degrees in the cone range.

Throughout the paper, $\pi:X\to\PP^1$ is a relatively minimal semistable
Jacobian elliptic surface over $\C$, with zero section $S=C_0$, fiber
class $F$, and
\[
 \chi=\chi(\OO_X)\ge3,\qquad H=S+\chi F.
\]
Then $q(X)=0$, $K_X\sim(\chi-2)F$, $S^2=-\chi$, and $S\cdot F=1$.
Let $s$ be the number of reducible fibers and $\rho(X)$ the Picard number.
We use $\Pic(X)$ for divisor classes. In the cases in which we compute
the full Cox ring, Lemma~\ref{lem:MW} gives a free Picard group. After
choosing divisor representatives for a free basis, its Cox ring is
\[
 \Cox(X)=\bigoplus_{D\in\Pic(X)}H^0(X,\OO_X(D)),
\]
with multiplication induced by these representatives. Its graded
isomorphism class is independent of the choices; see \cite{ADHL15,HK00}.
The section--fiber subring used below can be defined directly from $S,F$.
We call $D$ a \emph{necessary generator degree} if
$(R_+/R_+^2)_D\ne0$, where $R=\Cox(X)$ and $R_+$ is its ideal of nonzero
effective degrees; see Section~\ref{sec-nef-degree}.

If the Mordell--Weil group is finite and the reducible fibers have types
$I_{n_1},\ldots,I_{n_s}$, Shioda--Tate gives
\begin{equation}\label{eq:ST}
 \rho(X)=2+\sum_{i=1}^s(n_i-1).
\end{equation}
Put
\begin{equation}\label{eq:delta}
 N(\pi)=\sum_{i=1}^s n_i,\qquad
 \delta(\pi)=\sum_{i=1}^s\frac{\lfloor n_i^2/4\rfloor}{n_i},
\end{equation}
with both sums zero when $s=0$. Write the reducible fibers as
\[
 F_i=T_{i,0}+\cdots+T_{i,n_i-1},
\]
where $T_{i,0}$ meets $S$ and the components are cyclically ordered.
For $s\ge1$, let $\mathcal T$ be the set of all these components and put
\begin{equation}\label{eq:cone}
 \mathcal C_\pi=\R_{\ge0}[S]+\sum_{i,j}\R_{\ge0}[T_{i,j}].
\end{equation}
The following result is the cone theorem of \cite[Theorem 1.1]{LLL26}.
\begin{theorem}\label{cone-thm}
Suppose that $\MW(\pi)$ is finite and $s\ge1$.
\begin{enumerate}
\item One has $\NE(X)=\mathcal C_\pi$ if and only if
$\delta(\pi)\le\chi$. In particular, the Mori cone is rational
polyhedral in this range.
\item If $N(\pi)\le2\chi+3$, then $X$ is a Mori dream surface.
\end{enumerate}
\end{theorem}

We first state our finite criterion. Assume $s\ge1$, $\MW(\pi)$ finite,
and $\delta(\pi)\le\chi$. Lemma~\ref{lem:MW} then gives
$\MW(\pi)=0$ and an integral Picard basis. Define additive monoids
\[
\begin{split}
 M&=\{A\in\Pic(X): A\text{ is nef and }A\cdot S=0\},\\
 M^*&=\{0\}\cup\{A\in M:A\cdot F\ge2\},\\
 P&=M^*+\sum_{i=1}^s\sum_{j=1}^{n_i-1}\mathbb N[T_{i,j}],
\end{split}
\]
where $0\in\mathbb N$. A nonzero element of $P$ is
\emph{indecomposable} if it is not a sum of two nonzero elements of $P$.
Let $\mathcal G$ consist of its indecomposable elements of positive
fiber degree. This is a finite, numerically determined set. More
explicitly, put $\mathcal Q=\operatorname{Hilb}(M)\setminus\{H\}$.
Then $\mathcal G$ is obtained from
\[
 \mathcal Q\cup\{2H,3H\}\cup\{H+Q:Q\in\mathcal Q\}
\]
by retaining precisely the elements indecomposable in $P$.
These decompositions may involve nonidentity fiber components;
$\mathcal G$ is therefore distinct from the Hilbert basis of $M$.
The details are proved in Section~\ref{sec:semistable-case}.

\begin{definition}[Nonvanishing along the zero section]\label{def:nonvanishing}
Under the preceding hypotheses, we say that $X$ satisfies
\emph{nonvanishing along the zero section} if every $G\in\mathcal G$
admits a section $g_G\in H^0(X,\OO_X(G))$ with $g_G|_S\ne0$.
\end{definition}
Since $G\cdot S=0$ and $S\cong\PP^1$, this is equivalent to
surjectivity of $H^0(X,\OO_X(G))\to H^0(S,\OO_S(G))\cong\C$,
or to the existence of an effective divisor in $|G|$ disjoint from $S$.

\begin{theorem}\label{thm:criterion}
Assume that $s\ge1$, $\MW(\pi)$ is finite, and $\delta(\pi)\le\chi$.
The set $\mathcal G$ is finite, and every $G\in\mathcal G$ is nef,
with $G\cdot S=0$ and $G\cdot F\ge2$.
Suppose that $X$ satisfies nonvanishing along the zero section: for each $G\in\mathcal G$ there exists
\[
g_G\in H^0(X,\mathcal O_X(G))\qquad\text{with}\qquad g_G|_S\ne0.
\]
Then a minimal homogeneous generating set of $\Cox(X)$ consists of:
\begin{enumerate}
\item the defining sections $s_S$ and $s_{i,j}$ of $S$ and all the
components $T_{i,j}$, including $T_{i,0}$;
\item one chosen section $g_G$ for each $G\in\mathcal G$;
\item $\max\{0,2-s\}$ additional sections of degree $F$, chosen so that,
together with the products $p_i=\prod_{j=0}^{n_i-1}s_{i,j}$,
they span $H^0(X,\mathcal O_X(F))$.
\end{enumerate}
In particular, $\Cox(X)$ is finitely generated and $X$ is a Mori dream surface.
\end{theorem}

The proof uses the cone description, finite generation of the monoid $P$,
and induction by restriction to $S$. It does not assume the Mori dream
property or depend on a separate criterion for one reducible fiber.
The combinatorial construction of $\mathcal G$ and the geometric
nonvanishing condition are separate steps.

Our second result verifies the condition in a uniform numerical range.
The even-fiber case of the Mori dream conclusion was established in
\cite[Theorem 1.4]{LL26}; the case with odd reducible fibers was posed
there as part of Problem~1.9. The following theorem covers both cases
and also supplies minimal generators.
\begin{theorem}\label{thm:bound}
Assume that $\MW(\pi)$ is finite, $s\ge1$, and $N(\pi)\le2\chi+4$.
Then $X$ satisfies nonvanishing along the zero section, and $\Cox(X)$
has the finite minimal homogeneous generating set described in
Theorem~\ref{thm:criterion}. In particular, $X$ is a Mori dream surface,
and no nef and big integral class $D$ with $D\cdot S>0$ is a necessary
generator degree.
\end{theorem}
The cone hypothesis is automatic here, since
\[
 \delta(\pi)\le\frac{N(\pi)}4\le\frac\chi2+1<\chi.
\]
The numerical bound is sufficient for arbitrary semistable
configurations; we do not assert its optimality for every configuration.
For a single reducible fiber it gives a direct specialization of the
general criterion, recorded in Section~\ref{subsec:single-fiber}.

We next state the explicit presentations. The section--fiber subring
is governed by the Weierstrass equation and already gives the full
Cox ring in Picard number two.
\begin{theorem}
\label{thm-section-fiber}
Assume that $\MW(\pi)$ is finite. 
Choose short Weierstrass coordinates and homogeneous base coordinates $u,v$. Put
$
 A\in\C[u,v]_{4\chi}, B\in\C[u,v]_{6\chi}.
$
Then the following statements hold.
\begin{enumerate}
	\item[(a)] The subring of the Cox ring supported on $\Z S+\Z F$ is
\begin{equation}
\label{thm-eq:R0}\tag{\RzeroTag}
 R_0:=\bigoplus_{a,b\in\Z}H^0(X,\OO_X(aS+bF))
 \cong\frac{\C[u,v,z,x,y]}{(y^2-x^3-Axz^4-Bz^6)}
\end{equation}
where
\[
 \deg u=\deg v=F,\quad\deg z=S,\quad\deg x=2H,\quad\deg y=3H.
\]
For $a\geq0$, its component of degree $aS+bF$ has the direct-sum description
\begin{equation*}
	 \bigoplus_{\substack{i\geq0,\ j\in\{0,1\}\\2i+3j\leq a}}
 \C[u,v]_{b-\chi(2i+3j)}\,z^{a-2i-3j}x^iy^j.
\end{equation*}
For $a<0$ the component is zero. The five indicated generators of $R_0$ are minimal.
\item[(b)] Assume $\rho(X)=2$. Then \textup{(\RzeroTag)} is the full Cox ring. Its five indicated generators are minimal, and its unique minimal relation has degree $6H$. In particular, $X$ is a Mori dream surface.
\end{enumerate}
\end{theorem}

For configurations of $I_2$ fibers, simultaneous interpolation of the
nodal coordinates yields a complete-intersection presentation. The
interpolation range extends to $s=2\chi+1$, one fiber beyond the cone
range $s\le2\chi$; at that endpoint an additional negative bisection
appears.
 \begin{theorem}
\label{thm-I_2}
Assume that $\MW(\pi)$ is finite, and all $s\ge1$ reducible fibers have type $I_2$ and $s\le 2\chi+1$.
Let $\ell_1,\ldots,\ell_s$ be pairwise nonproportional linear forms defining their base points, and let $E_i$ be the nonidentity component of the $i$th fiber.
Put $\mathcal E=\prod_i e_i$ and $\mathcal T=\prod_i t_i$. Then the following statements hold.
\begin{enumerate}
	\item 
\begin{equation*}
 \Cox(X)\cong
 \frac{\C[u,v,z,e_1,t_1,\ldots,e_s,t_s,\xi,\eta]}
 {(\ell_i(u,v)-e_it_i\ (1\leq i\leq s),\ G_s)},
\end{equation*}
where
\[
 G_s=\eta^2-\mathcal E\xi^3-3r\xi^2z^2
            -\mathcal T\alpha\xi z^4-\mathcal T^2\beta z^6,
\]
and $r,\alpha,\beta$ are the polynomials in $u,v$ of \eqref{eq:multiab}. The degrees are
\begin{gather*}
 \deg u=\deg v=F,\quad \deg z=S,\quad
 \deg e_i=E_i,\quad \deg t_i=F-E_i,\\
 \deg\xi=2H-\sum_iE_i,\qquad
 \deg\eta=3H-\sum_iE_i.
\end{gather*}
For $s=1$, eliminating one base variable gives six minimal generators and one relation. For $s\ge2$, eliminating both base variables gives $2s+3$ minimal generators and $s-1$ minimal relations. In particular, $X$ is a Mori dream surface.
\item If $s=2\chi+1$, the divisor of $\xi$ is an irreducible horizontal bisection $C$ with $C^2=-2$ and $p_a(C)=\chi-2$.
In particular, the visible fiber components and the zero section no longer generate the effective cone even though the Cox ring is finitely generated.
\item Put $D_\xi=2H-\sum_iE_i$. $\Eff(X)=\mathrm{cone}(S, E_i, F-E_i, D_\xi)$.
For $s\leq2\chi$, the last class is redundant. For $s=2\chi+1$, a class $D=aS+bF-\sum_i d_iE_i$ is nef precisely when
\[
 0\leq d_i\leq a/2\quad\text{for every }i,\qquad
 b\geq\chi a,\qquad b\geq\sum_i d_i.
\]
\end{enumerate}
\end{theorem}

The bound $s\le2\chi+1$ is an interpolation bound, conditional on the
existence of a surface with the stated configuration; no existence
claim is made for every pair $(\chi,s)$. Section~\ref{subsec:multigraded}
places these presentations in a common framework of simultaneous
exceptional valuations.

The exclusion of ample generator degrees holds throughout the cone
range, independently of the nonvanishing condition and of finite
generation.
\begin{theorem}
\label{thm:no-ample}
Assume that $\MW(\pi)$ is finite and $\delta(\pi)\leq\chi$. 
For every ample integral divisor $D$ on $X$, the following statements hold.
\begin{enumerate}
\item[(i)] One has $H^1(X,\OO_X(D-2F))=0$, and the multiplication map
\[
 H^0(X,\OO_X(D-F))\otimes H^0(X,\OO_X(F))
 \longrightarrow H^0(X,\OO_X(D))
\]
is surjective.
\item[(ii)] Put $d=D\cdot S$ and $D_0=D-dF$. Then $D_0$ is nef, $D_0\cdot S=0$, and multiplication induces a surjection
\[
 H^0(X,\OO_X(D_0))\otimes
 \operatorname{Sym}^{d}H^0(X,\OO_X(F))
 \longrightarrow H^0(X,\OO_X(D)).
\]
\item[(iii)] The Cox ring has no generator in any ample degree. If it is finitely generated, every minimal homogeneous generator is either, up to a nonzero scalar, the defining section of an irreducible negative curve or has degree on the boundary of $\Nef(X)$.
\end{enumerate}
\end{theorem}

The cone condition is the input to the vanishing argument. Examples in
\cite[Theorem 1.5]{LLL26} and \cite[Theorem 1.7]{LL26} show that surfaces
in this range can fail to be Mori dream surfaces. Conversely,
Corollary~\ref{cor:explicit-no-ample} excludes ample degrees at the
$I_2$ endpoint $s=2\chi+1$, beyond this range, even though the uniform
vanishing fails there (Example~\ref{ex:endpoint-vanishing}).

These results leave three related questions. First, reduction to the
nef boundary does not identify all necessary degrees: a decomposition
of a divisor class need not give a surjective multiplication map of
sections. Even on $S^\perp$, the integral nef Hilbert basis need not
contain the generator degrees, as $2H,3H$ in Picard number two show.
\begin{problem}\label{prop-main}
Let $X$ be a Mori dream surface in the above semistable Jacobian setting,
with finite Mordell--Weil group.
\begin{enumerate}
\item If $\delta(\pi)>\chi$, determine whether ample necessary generator
degrees can occur, and give conditions excluding them.
\item Determine the necessary non-ample nef degrees on the boundary
faces of $\Nef(X)$.
\end{enumerate}
\end{problem}

In all the explicit presentations and in Theorem~\ref{thm:criterion},
every nef and big generator degree is orthogonal to $S$. This suggests
a more specific question.
\begin{question}\label{ques:positive-section}
Does there exist a Mori dream surface in the above semistable Jacobian
setting, with finite Mordell--Weil group, and a nef and big necessary
generator degree $D$ with $D\cdot S>0$? Can such an example occur when
$\delta(\pi)\le\chi$?
\end{question}
Theorem~\ref{thm:bound} excludes such examples when
$N(\pi)\le2\chi+4$. Finally, the finite criterion separates the
classification problem from the construction of candidate degrees.
\begin{problem}\label{prob:nonvanishing}
Among surfaces with $s\ge1$, finite Mordell--Weil group, and
$\delta(\pi)\le\chi$, characterize those satisfying nonvanishing along
the zero section. Determine how this condition depends on the fiber
configuration and on the surface within that configuration. Is it
also necessary for finite generation of the Cox ring?
\end{problem}

Section~\ref{Sect-section-fiber} establishes the Picard-lattice and
cohomological preliminaries and computes the section--fiber subring,
including the Picard-number-two case. Section~\ref{Sect-I_2} gives
the complete-intersection presentations for $I_2$ configurations,
their cones and cohomology, and the description by exceptional
valuations. Section~\ref{sec:ample} proves the uniform vanishing
theorem and explains its limits outside the cone range.
Section~\ref{sec:semistable-case} proves the finite criterion and
the nonvanishing theorem for $N(\pi)\le2\chi+4$, then specializes
the criterion to a single reducible fiber.

\smallskip\noindent\textbf{Acknowledgments.}
Antonio Laface was partially supported by Proyecto FONDECYT Regular
n.~1230287. Sichen Li thanks Rong Du for constant encouragement.

\smallskip\noindent\textbf{AI disclosure.}
During the preparation of this work, the authors used ChatGPT (OpenAI)
to assist with editing and revising the manuscript, discussing
mathematical arguments, and working out and checking examples. The
authors take full responsibility for the content of the paper.

\section{The section--fiber subring}
\label{Sect-section-fiber}
\subsection{Vanishing of the Mordell-Weil group}

\begin{lemma}\label{lem:MW}
If $\MW(\pi)$ is finite and $\delta(\pi)<2\chi$, then $\MW(\pi)=0$. In that case $S,F$, and all nonidentity fiber components form an integral basis of $\Pic(X)$.
\end{lemma}
\begin{proof}
For a nonzero torsion section $P$, the height formula gives
\[
 0=\langle P,P\rangle=2\chi+2P\cdot S-\sum_i\operatorname{contr}_i(P).
\]
For a fiber of type $I_n$ each local contribution is at most $\lfloor n^2/4\rfloor/n$. Since $P\cdot S\geq0$, this contradicts $\delta(\pi)<2\chi$. The integral-basis assertion follows from the exact description of the N\'eron--Severi group modulo the trivial lattice by the Mordell--Weil group, and $q(X)=0$. See \cite[Theorem 6.3, Corollary 6.13, Section 11]{SS10} and \cite[Lemma 2.1]{LL26}.
\end{proof}
\subsection{Homogeneous generators and multiplication maps}
\label{sec-nef-degree}
For a surface with free Picard group write $R=\Cox(X)$ and
\[
 R_+=\bigoplus_{\substack{D\ne0\\H^0(X,\OO_X(D))\ne0}} R_D.
\]
Intersection with an ample divisor gives a positive integral grading on every nonzero effective degree. When $R$ is finitely generated, the degree multiplicities in a minimal homogeneous generating set are therefore the dimensions of the homogeneous components of $R_+/R_+^2$. We say that $R$ has a generator in degree $D$ when $(R_+/R_+^2)_D\ne0$.

If homogeneous elements $f_1,\ldots,f_N$ generate $R$, every effective class is a nonnegative integral combination of their degrees. Indeed, a nonzero homogeneous section is a polynomial in the $f_i$, and some monomial in its homogeneous expression is nonzero. In particular,
\[
 \Eff(X)=\operatorname{cone}(\deg f_1,\ldots,\deg f_N).
\]
See also \cite[Proposition 1.3]{ADL21}.

\begin{proposition}\label{prop:koszul}
Let $B_1,B_2$ be nonzero effective divisors with disjoint supports, and let $s_1,s_2$ be their defining sections. For an effective class $D$ with $D\not\sim B_1,B_2$, consider
\[
 H^0(X,\OO_X(D-B_1))\oplus H^0(X,\OO_X(D-B_2))
 \xrightarrow{\ (s,t)\mapsto s_1s+s_2t\ }
 H^0(X,\OO_X(D)).
\]
If this map is surjective, there is no generator in degree $D$. In particular, it is enough that $H^1(X,\OO_X(D-B_1-B_2))=0$.
\end{proposition}
\begin{proof}
Disjointness makes the twisted Koszul sequence exact:
\[
 0\longrightarrow\OO_X(D-B_1-B_2)
 \longrightarrow\OO_X(D-B_1)\oplus\OO_X(D-B_2)
 \longrightarrow\OO_X(D)\longrightarrow0.
\]
The long exact cohomology sequence gives the sufficient condition for surjectivity. Every nonzero product in the image has two factors of nonzero effective degree, so the image lies in $(R_+^2)_D$. This proves the assertion; compare \cite[Lemma 1.1]{AGL16}.
\end{proof}

\subsection{Fixed negative components}
\begin{proposition}\label{prop:negative-nef}
Let $Y$ be a smooth projective surface with free Picard group. Every necessary generator degree of $\Cox(Y)$ is either the class of an irreducible negative curve, in which case its section space is spanned by the defining section of that curve, or is nef.
\end{proposition}
\begin{proof}
If an effective divisor class $D$ is not nef, choose an irreducible curve $C$ with $D\cdot C<0$. Every effective representative of $D$ contains $C$, and $C^2<0$. Multiplication by its defining section induces
\[
 H^0(Y,\OO_Y(D-C))\xrightarrow{\ \cdot s_C\ }H^0(Y,\OO_Y(D)).
\]
This is an isomorphism. If $D\ne C$, both factors have nonzero effective degree, so a section of degree $D$ is decomposable in the positive-degree ideal. If $D=C$, the space of sections is one-dimensional, generated by $s_C$. Conversely $s_C$ cannot be a product of nonconstant homogeneous sections, since its zero divisor is irreducible.
\end{proof}

\begin{remark}
This statement concerns nef \emph{degrees}, not merely the Hilbert basis of the integral nef monoid. A decomposition $D=D_1+D_2$ does not imply surjectivity of $H^0(D_1)\otimes H^0(D_2)\to H^0(D)$. Both explicit cases below exhibit indispensable generators in decomposable nef degrees.
\end{remark}
\subsection{Vertical classes and cohomology in the section--fiber lattice}

\begin{proposition}\label{prop:vertical-nef}
Let $D$ be a nonzero effective nef divisor on $X$. If $D\cdot F=0$,
then $D\sim bF$ for an integer $b>0$. Moreover,
\[
 \bigoplus_{b\geq0}H^0(X,\OO_X(bF))\cong\C[u,v],
 \qquad \deg u=\deg v=F.
\]
\end{proposition}
\begin{proof}
Every component of $D$ is vertical, since $F$ is nef and has positive
intersection with every horizontal curve. On a reducible fiber
$F_p=\sum_jC_j$, nefness and $D\cdot F_p=0$ imply $D\cdot C_j=0$
for every $j$. The kernel of the intersection matrix of the components
of a connected fiber is generated by its multiplicity vector. The fibers
are reduced, so the part of $D$ supported on $F_p$ is an integral
multiple of $F_p$. The same conclusion is immediate on an irreducible
fiber. Thus $D$ is a sum of full fibers, proving the first assertion.
The second follows from $\pi_*\OO_X=\OO_{\PP^1}$, the projection
formula, and the usual section ring of $\OO_{\PP^1}(1)$.
\end{proof}

\begin{remark}
Nefness is essential in the first assertion: an individual component
of a reducible fiber has intersection zero with $F$ but is not linearly
equivalent to a multiple of $F$.
\end{remark}

Write $z$ for the canonical section of $\OO_X(S)$.

\begin{lemma}\label{lem:pole-filtration}
For $a\geq1$ one has
\[
 R^1\pi_*\OO_X(aS)=0,\qquad
 \pi_*\OO_X(aS)\cong
 \OO_{\PP^1}\oplus\bigoplus_{n=2}^{a}\OO_{\PP^1}(-n\chi).
\]
The sum is empty when $a=1$. With compatible Weierstrass coordinates,
the summand of pole order $n=2i+3j$, with $j\in\{0,1\}$, is locally represented
by $z^{a-2i-3j}x^iy^j$; the summand of pole order zero is represented
by $z^a$.
\end{lemma}
\begin{proof}
Since $\omega_T\cong\OO_T$, on every fiber $T$, duality identifies $H^1(T,\OO_T(aS))$ with the
dual of $H^0(T,\OO_T(-aS))$.
The latter space is zero. For a smooth or integral nodal fiber this follows from
negative degree. On an $I_n$ fiber the section has negative degree on
the identity component, so vanishes there; it has degree zero on all
other components, and the vanishing propagates across the connected
cycle. 
Cohomology and base change therefore give the asserted higher
direct-image vanishing. 
The canonical section of $S$ spans the one-dimensional space $H^0(T,\OO_T(S))$ on every fiber, so
$\pi_*\OO_X(S)=\OO_{\PP^1}$.

Set $V_a=\pi_*\OO_X(aS)$. Restriction to $S$ gives, for $a\geq2$,
\[
 0\longrightarrow V_{a-1}\longrightarrow V_a
 \longrightarrow\OO_{\PP^1}(-a\chi)\longrightarrow0.
\]
Inductively the extension splits, because its obstruction group is a
sum of groups $$H^1(\PP^1,\OO((a-n)\chi)),$$
 where
$n\in\{0,2,\ldots,a-1\}$, and these groups vanish. This proves the
splitting with the indicated negative signs.

Choose the pole-order-two and pole-order-three sections $x$ and $y$,
twisted by $2\chi F$ and $3\chi F$, so that their leading coefficients
are nonzero constants. For every $n\geq2$ there is a unique expression
$n=2i+3j$ with $i\geq0$ and $j\in\{0,1\}$. The corresponding monomial
has a nonzero leading coefficient in the quotient of pole order $n$.
Consequently these monomials give a basis compatible with the
filtration, as claimed. Completing the square and cube give the short
Weierstrass relation in characteristic zero without changing the
leading pole orders.
\end{proof}

\begin{corollary}\label{cor:section-fiber-cohomology}
For any $b\in\Z$ and $a\geq1$,
\begin{align*}
 h^1(X,\OO_X(aS+bF))
 &=\max\{-b-1,0\}+
   \sum_{n=2}^{a}\max\{n\chi-b-1,0\},\\
 h^2(X,\OO_X(aS+bF))&=0.
\end{align*}
Thus $h^1(aS+bF)=0$ precisely when $b\geq-1$ for $a=1$, and
precisely when $b\geq a\chi-1$ for $a\geq2$. In particular it
vanishes whenever $aS+bF$ is nef and $a>0$. By contrast,
\[
 h^1(X,\OO_X(bF))=
 \max\{-b-1,0\}+\max\{b-\chi+1,0\}.
\]
\end{corollary}
\begin{proof}
Apply the Leray spectral sequence, the projection formula, and Lemma~\ref{lem:pole-filtration}.
The nef assertion follows from $(aS+bF)\cdot S=b-a\chi\geq0$.
For $a=0$, use $\pi_*\OO_X=\OO_{\PP^1}$ and
$R^1\pi_*\OO_X=\OO_{\PP^1}(-\chi)$ in the Leray sequence.
\end{proof}
\subsection{The section--fiber subring}
Let $f:X\to W$ contract the fiber components disjoint from $S$. The surface $W$ is the normal Weierstrass model, and $S$ descends to a Cartier divisor, still denoted $S$. Choose short Weierstrass coordinates and homogeneous base coordinates $u,v$. Put
\[
 A\in\C[u,v]_{4\chi},\qquad B\in\C[u,v]_{6\chi}.
\]
Here $\C[u,v]_d=0$ for $d<0$.

\begin{proposition}\label{prop:subring}
The subring of the Cox ring supported on $\Z S+\Z F$ is
\begin{equation}\label{eq:R0}
 R_0:=\bigoplus_{a,b\in\Z}H^0(X,\OO_X(aS+bF))
 \cong\frac{\C[u,v,z,x,y]}{(y^2-x^3-Axz^4-Bz^6)},
\end{equation}
where
\[
 \deg u=\deg v=F,\quad\deg z=S,\quad\deg x=2H,\quad\deg y=3H.
\]
For $a\geq0$, its component of degree $aS+bF$ has the direct-sum description
\begin{equation}\label{eq:normalforms}
 \bigoplus_{\substack{i\geq0,\ j\in\{0,1\}\\2i+3j\leq a}}
 \C[u,v]_{b-\chi(2i+3j)}\,z^{a-2i-3j}x^iy^j.
\end{equation}
For $a<0$ the component is zero.
\end{proposition}
\begin{proof}
Since $aS+bF$ is pulled back from a Cartier divisor on $W$ and $f_*\OO_X=\OO_W$, its sections can be computed on $W$. Lemma~\ref{lem:pole-filtration} identifies the pole-order basis and gives
\begin{equation}\label{eq:pushforward}
 \pi_*\OO_X(aS)=\OO\oplus\bigoplus_{n=2}^{a}\OO(-n\chi)
 \qquad(a\geq1),
\end{equation}
while $\pi_*\OO_X=\OO$. See also \cite[Lemma 3.5]{LMS09}. Twisting by $\OO_{\PP^1}(b)$ proves \eqref{eq:normalforms}. In particular, the direct-image summands have negative degrees.

The defining equation is monic of degree two in $y$. 
Consequently its polynomial quotient has unique normal forms with $y$-exponent at most one, whose homogeneous components are exactly \eqref{eq:normalforms}. The natural graded map is therefore an isomorphism in every degree. For $a<0$, the restriction of a putative section to the generic fiber would have negative degree and must vanish. This completes the isomorphism proof. The defining relation has no ordinary linear term, so all five generators are minimal.
\end{proof}
\subsection{Picard number two}
\begin{theorem}
\label{thm:rho2}
Assume $\rho(X)=2$. Then $\Pic(X)=\Z S\oplus\Z F$, and \eqref{eq:R0} is the full Cox ring. Its five indicated generators are minimal, and its unique minimal relation has degree $6H$. In particular, $X$ is a Mori dream surface.
\end{theorem}
\begin{proof}
Shioda--Tate, including the Mordell--Weil rank term, implies that there are no reducible fibers and the Mordell--Weil group has rank zero. Its torsion is zero by Lemma~\ref{lem:MW}, with $\delta=0$. Thus $S,F$ form an integral Picard basis and Proposition~\ref{prop:subring} computes the full ring.

The space $H^0(F)$ has dimension two, so two generators in degree $F$ are necessary. The defining section $z$ of $S$ is necessary. In degree $2H$, the products in $u,v,z$ span precisely
\[
 z^2\C[u,v]_{2\chi},
\]
and $x$ contributes one additional direction. In degree $3H$, the products in $u,v,z,x$ span
\[
 z^3\C[u,v]_{3\chi}\oplus xz\C[u,v]_{\chi},
\]
and $y$ contributes one more direction. 
Thus all five generators are minimal. 
The normal-form proof establishes that the displayed equation generates the entire kernel. 
Thus, $R$ is finitely generated, and $X$ is a Mori dream surface.
\end{proof}

\begin{corollary}\label{cor:rho2cones}
Under the hypotheses of Theorem~\ref{thm:rho2},
\[
 \Eff(X)=\R_{\geq0}S+\R_{\geq0}F,\qquad
 \Nef(X)=\R_{\geq0}F+\R_{\geq0}H.
\]
The integral nef monoid has Hilbert basis $\{F,H\}$. The complete linear system $|H|$ has fixed component $S$, whereas $|mH|$ is base-point-free for every $m\geq2$.
\end{corollary}
\begin{proof}
All generator degrees belong to the cone spanned by the effective classes $S,F$, so this is the effective cone. Its dual is described by $a\geq0$ and $b\geq\chi a$ for $aS+bF$, giving the asserted nef cone and integral monoid. Formula~\eqref{eq:normalforms} gives $H^0(H)=z\C[u,v]_\chi$. Off $S$, the sections $z^m\C[u,v]_{m\chi}$ have no common zero. Along $S$, the sections $x$ and $y$ do not vanish; for even $m\geq2$ use $x^{m/2}$, and for odd $m\geq3$ use $x^{(m-3)/2}y$.
\end{proof}

\begin{corollary}\label{cor:rho2hilbert}
Writing $s,t$ for the formal degrees $S,F$, the Hilbert series is
\[
 \frac{1-s^6t^{6\chi}}
 {(1-t)^2(1-s)(1-s^2t^{2\chi})(1-s^3t^{3\chi})}.
\]
For $a\geq0$,
\begin{equation}\label{eq:h0rho2}
 h^0(aS+bF)=\pos{b+1}+\sum_{n=2}^{a}\pos{b-n\chi+1}.
\end{equation}
For $a\geq1$,
\[
 h^1(aS+bF)=\pos{-b-1}+\sum_{n=2}^{a}\pos{n\chi-b-1},
 \qquad h^2(aS+bF)=0.
\]
\end{corollary}
\begin{proof}
The Hilbert series follows from the homogeneous hypersurface presentation, and \eqref{eq:h0rho2} from the normal forms. For $a\geq1$, cohomology and base change give $R^1\pi_*\OO_X(aS)=0$, because the fibers are integral of arithmetic genus one. Apply the Leray spectral sequence and \eqref{eq:pushforward} to obtain the desired result.
\end{proof}

\begin{remark}\label{rem:false}
The only negative irreducible curve when $\rho=2$ is $S$. Indeed, another horizontal irreducible curve $C=aS+bF$ satisfies $a>0$ and $b\geq\chi a$, so $C^2\geq\chi a^2>0$; vertical irreducible curves are fibers. If $\chi\geq4$, sections of degrees $S$ and $K_X=(\chi-2)F$ cannot even generate the nonzero component of degree $F$. For $\chi=3$, the degree $F$ is available, but the generators $x$ and $y$ of degrees $2H,3H$ are still missing. Thus the defining sections of negative curves, together with sections of degree $K_X$ or $F$, do not suffice to generate the Cox ring. Likewise the nef Hilbert basis $\{F,H\}$ does not suffice.
\end{remark}
The following result answers \cite[Question 5.10]{HL27} for Mori dream Jacobian elliptic surfaces of Picard number two.
\begin{corollary}\label{cor:bcp}
If $\rho(X)=2$, then $H^1(X,\OO_X(C))=0$ for every reduced irreducible curve $C$.
\end{corollary}
\begin{proof}
For $C=S$ use Corollary~\ref{cor:rho2hilbert}. 
For another horizontal curve, $C=aS+bF$ with $a\geq1$ and $b\geq\chi a$, and the same formula gives zero. A vertical irreducible curve has class $F$.
By the Leray spectral sequence, $\pi_*\OO_X=\OO$, and $R^1\pi_*\OO_X=\OO(-\chi)$ give $h^1(F)=h^1(\OO_{\PP^1}(1))+h^0(\OO_{\PP^1}(1-\chi))=0$.
\end{proof}
\begin{proof}[Proof of Theorem \ref{thm-section-fiber}]
It follows from Proposition \ref{prop:subring} and Theorem \ref{thm:rho2}.
\end{proof}
\section{Surfaces with reducible fibers of type \texorpdfstring{$I_2$}{I2}}
\label{Sect-I_2}
\subsection{Picard number three: one fiber of type \texorpdfstring{$I_2$}{I2}}
Assume that $\rho(X)=3$ and $\MW(\pi)$ is finite. By \eqref{eq:ST}, exactly one fiber is reducible and it has type $I_2$. Write it as $E_0+E$, with $S\cdot E_0=1$ and $S\cdot E=0$. Lemma~\ref{lem:MW} gives
\[
 \Pic(X)=\Z S\oplus\Z F\oplus\Z E,
 \qquad E^2=-2,\quad E\cdot F=E\cdot S=0.
\]
Here $E_0=F-E$. The Weierstrass model has a unique singular point, an ordinary double point obtained by contracting $E$.
\par\smallskip\noindent\textbf{Normalized Weierstrass coefficients.}
Choose the base coordinate so that the $I_2$ fiber lies over $u=0$ and its nodal point, on $v=z=1$, is $(x,y)=(r_0,0)$, with $r_0\ne0$. Set
\[
 r=r_0v^{2\chi},\qquad X_1=x-rz^2.
\]
There are homogeneous polynomials
\begin{equation}\label{eq:alphabeta}
 \alpha=\frac{A+3r^2}{u}\in\C[u,v]_{4\chi-1},
 \qquad
 \beta=\frac{B+rA+r^3}{u^2}\in\C[u,v]_{6\chi-2}.
\end{equation}
To justify the second divisibility, the double root gives $A(0,1)=-3r_0^2$ and $B(0,1)=2r_0^3$. The vanishing of the first derivative of $4A^3+27B^2$ gives $(B+rA+r^3)'(0,1)=0$. In these coordinates the equation is
\begin{equation}\label{eq:shifted}
 y^2=X_1^3+3rX_1^2z^2+u\alpha X_1z^4+u^2\beta z^6.
\end{equation}
Put $\alpha_0=\alpha(0,1)$ and $\beta_0=\beta(0,1)$. The exact discriminant expansion is
\begin{align}\label{eq:disc}
 4A^3+27B^2
 &=9r^2(12r\beta-\alpha^2)u^2
 +(4\alpha^3-54r\alpha\beta)u^3+27\beta^2u^4.
\end{align}
The fiber has type exactly $I_2$, hence
\begin{equation}\label{eq:nondeg}
 \alpha_0^2-12r_0\beta_0\ne0.
\end{equation}
This is also precisely the nondegeneracy condition for the quadratic tangent cone of the surface singularity. See \cite[Sections 4.2--4.3]{SS10} for the local resolution of multiplicative fibers.

\begin{theorem}\label{thm:rho3}
Under these assumptions the Cox ring is
\begin{equation}\label{eq:rho3ring}
 \Cox(X)\cong\C[v,z,e,w,\xi,\eta]/(G),
\end{equation}
where
\begin{align}\label{eq:G}
 G={}&\eta^2-e\xi^3-3r\xi^2z^2
       -w\alpha(ew,v)\xi z^4-w^2\beta(ew,v)z^6,
 \qquad r=r_0v^{2\chi},
\end{align}
and the degrees are
\begin{center}
\begin{tabular}{c|cccccc}
variable & $z$ & $e$ & $w$ & $v$ & $\xi$ & $\eta$\\ \hline
degree & $S$ & $E$ & $F-E$ & $F$ & $2H-E$ & $3H-E$
\end{tabular}
\end{center}
These six generators are minimal; the unique minimal relation has degree $6H-2E$. In particular $X$ is a Mori dream surface.
\end{theorem}

\begin{proof}
We prove that the coordinate substitution computes the \emph{whole} Cox ring. Let $R_0$ be the ring of Proposition~\ref{prop:subring}, written with $X_1$ in place of $x$, and set
\[
 I=(u,X_1,y)\subset R_0.
\]
This is the homogeneous ideal of the node in the chart where $v,z$ are invertible. The $I$-adic associated graded ring is
\begin{equation}\label{eq:gradedI}
 \gr_I R_0\cong\frac{\C[v,z][U,T,Y]}{(q)},
\end{equation}
where
\[
 q=Y^2-3r_0v^{2\chi}T^2z^2
 -\alpha_0v^{4\chi-1}UTz^4
 -\beta_0v^{6\chi-2}U^2z^6.
\]
Indeed, in the ambient polynomial ring with ideal $(u,X_1,y)$, the initial form of the relation \eqref{eq:shifted} is $q$. The ambient associated graded ring is a domain, so the initial form of any multiple of the relation is the product of the initial forms. Thus the initial ideal is generated by $q$.

Over $\C(v,z)$, the quadratic form $q$ in $U,T,Y$ has rank three by \eqref{eq:nondeg}; it is therefore irreducible. Since it is monic in $Y$, Gauss's lemma implies irreducibility over $\C[v,z]$, and \eqref{eq:gradedI} is a domain. Consequently $I$-adic order is additive on nonzero products and is unchanged by localizing at $v$ and $z$. On that chart, after trivializing the two divisor degrees, \eqref{eq:shifted} is the local equation of an ordinary double point. Blowing up its maximal ideal resolves it, with irreducible exceptional conic $E$. The order along this exceptional divisor equals the maximal-ideal order, since the associated graded ring is the domain defining that conic. Hence, for a homogeneous section $s\in R_0$,
\begin{equation}\label{eq:order}
 \ord_E(s)\geq m\quad\Longleftrightarrow\quad s\in I^m
 \qquad(m\geq0).
\end{equation}
There is no additional saturation along $v=0$ or $z=0$: their degree-zero initial forms are nonzero in the domain \eqref{eq:gradedI}.

For $D=aS+bF$, equation \eqref{eq:order} identifies
\[
 H^0(X,\OO_X(D-mE))=(I^m)_D.
\]
For $c>0$, the intersection $(D+cE)\cdot E=-2c$ forces $E$ to be a fixed component. Removing it repeatedly gives
\[
 H^0(X,\OO_X(D+cE))=e^cH^0(X,\OO_X(D)),
\]
where $e$ is the defining section of $E$. As $S,F,E$ form an integral Picard basis, these identifications give the extended Rees algebra
\begin{equation}\label{eq:Rees}
 \Cox(X)=R_0[e,Ie^{-1}].
\end{equation}
Its generators may be written
\[
 w=u/e,\qquad \xi=X_1/e,\qquad \eta=y/e,
\]
together with $e,v,z$. Substituting $u=ew$, $X_1=e\xi$, $y=e\eta$ into \eqref{eq:shifted} and dividing by $e^2$ gives $G$.

For completeness, this substitution introduces no extra relation. After inverting $e$, the quotient by $G$ is $R_0[e,e^{-1}]$. Moreover $G$ is not divisible by $e$ in the ambient polynomial ring, so $(G):e^\infty=(G)$. The quotient embeds in its localization at $e$ and hence in $R_0[e,e^{-1}]$; its image is exactly \eqref{eq:Rees}. This proves \eqref{eq:rho3ring}.

All six variable degrees are nonzero effective classes, and intersection with an ample divisor gives a positive integer grading. At the homogeneous origin, $G$ has no ordinary linear term, so the cotangent space of the augmentation ideal has dimension six. Thus no generator can be removed. Homogeneity gives relation degree $6H-2E$, and the kernel is the principal ideal already determined.
\end{proof}

\begin{corollary}\label{cor:rho3cones}
The effective and nef cones are
\[
 \Eff(X)=\operatorname{cone}(S,E,F-E),\qquad
 \Nef(X)=\operatorname{cone}(F,H,Q),\quad Q=2H-E.
\]
The integral nef monoid is freely generated by $F,H,Q$.
\end{corollary}
\begin{proof}
All six degrees in Theorem~\ref{thm:rho3} lie in the cone spanned by $S,E,F-E$, and these three effective classes occur as generator degrees. Duality gives, for $D=aS+bF+cE$,
\[
 D\cdot S=b-\chi a\geq0,\quad
 D\cdot E=-2c\geq0,\quad D\cdot(F-E)=a+2c\geq0.
\]
Writing $d=-c$, the integral solutions are uniquely
\[
 D=(a-2d)H+dQ+(b-\chi a)F,
 \quad d\geq0,\quad a\geq2d,\quad b\geq\chi a.
\]
\end{proof}

\begin{corollary}\label{cor:rho3h0}
For $a\geq0$ and arbitrary integers $b,c$,
\begin{equation}\label{eq:h0rho3}
 h^0(aS+bF+cE)=
 \sum_{\substack{i\geq0,\ j\in\{0,1\}\\2i+3j\leq a}}
 \pos{b-\chi(2i+3j)-\max\{0,-c-i-j\}+1}.
\end{equation}
For $a<0$ the dimension is zero. If $D=aS+bF-dE$ is nef and $a>0$, then
\begin{align*}
 h^0(D)&=a(b+1)-\frac{\chi a(a+1)}2+\chi-d^2,\\
 h^1(D)&=h^2(D)=0.
\end{align*}
\end{corollary}
\begin{proof}
Since $G$ is monic in $\eta$, a vector-space basis is given by monomials
\[
 e^p w^q v^r z^k\xi^i\eta^j,\qquad j=0,1.
\]
In degree $aS+bF+cE$, one has
\[
 k=a-2i-3j,\quad q+r=b-\chi(2i+3j),\quad p=c+q+i+j.
\]
The possible $q$ are the integers from $\max\{0,-c-i-j\}$ to $b-\chi(2i+3j)$, proving the formula.

For a nef divisor, $a\geq2d$ and $b\geq\chi a$. Relative to the rank-two sum, the total correction is
\[
 d+2\sum_{r=1}^{d-1}(d-r)=d^2.
\]
Summing the remaining terms gives the displayed expression. Since
\[
 D^2=-\chi a^2+2ab-2d^2,\qquad K_X\cdot D=(\chi-2)a,
\]
this expression equals $\chi(\OO_X(D))$. Finally $h^2(D)=0$, because $(K_X-D)\cdot F=-a<0$ and $F$ is nef. Riemann--Roch then gives $h^1(D)=0$.
\end{proof}

\begin{remark}
For $\chi=3$, the generator degrees in the basis $(S,F,E)$ are
\[
 (1,0,0),\ (0,0,1),\ (0,1,-1),\ (0,1,0),\ (2,6,-1),\ (3,9,-1),
\]
and the relation degree is $(6,18,-2)$. For example $h^0(H)=4$, $h^0(Q)=7$, and $h^0(H+Q)=14$. The degree of $\eta$ is the decomposable nef class $H+Q$, although $\eta$ is indispensable. These are universal calculations under the theorem's hypotheses, not an assertion that an arbitrary choice of coefficients has Picard number three.
\end{remark}
\subsection{Several fibers of type \texorpdfstring{$I_2$}{I2}}
Assume that all $s\geq1$ reducible fibers have type $I_2$, that $\MW(\pi)$ is finite, and that $s\leq2\chi+1$. Let $\ell_1,\ldots,\ell_s$ be pairwise nonproportional linear forms defining their base points, and let $E_i$ be the component of the $i$th fiber disjoint from $S$. Since $\delta=s/2<2\chi$, Lemma~\ref{lem:MW} gives trivial Mordell--Weil group and the integral basis $S,F,E_1,\ldots,E_s$.
\par\smallskip\noindent\textbf{Simultaneous translation of the nodes.}
The double root of the Weierstrass cubic at each reducible fiber is a prescribed nonzero value in the fiber of $\OO_{\PP^1}(2\chi)$. Interpolation gives a section
\[
 r\in\C[u,v]_{2\chi}
\]
having all these values: the evaluation map at $s$ distinct points is surjective because $H^1(\PP^1,\OO(2\chi-s))=0$. Put $L=\prod_i\ell_i$. As in \eqref{eq:alphabeta}, now at every nodal point, we obtain
\begin{equation}\label{eq:multiab}
 A+3r^2=L\alpha,\qquad B+rA+r^3=L^2\beta,
 \quad
 \deg\alpha=4\chi-s,\quad \deg\beta=6\chi-2s.
\end{equation}
The derivative of $r$ does not affect the second-order divisibility, since its coefficient is $A+3r^2$, which vanishes at each node. For $X_1=x-rz^2$, the shifted relation is
\begin{equation}\label{eq:multishift}
 y^2=X_1^3+3rX_1^2z^2+L\alpha X_1z^4+L^2\beta z^6.
\end{equation}

\begin{theorem}\label{thm:multiI2}
Under these assumptions, put $\mathcal E=\prod_i e_i$ and $\mathcal T=\prod_i t_i$. Then
\begin{equation}\label{eq:multiI2}
 \Cox(X)\cong
 \frac{\C[u,v,z,e_1,t_1,\ldots,e_s,t_s,\xi,\eta]}
 {(\ell_i(u,v)-e_it_i\ (1\leq i\leq s),\ G_s)},
\end{equation}
where
\[
 G_s=\eta^2-\mathcal E\xi^3-3r\xi^2z^2
            -\mathcal T\alpha\xi z^4-\mathcal T^2\beta z^6,
\]
and $r,\alpha,\beta$ are the polynomials in $u,v$ of \eqref{eq:multiab}. The degrees are
\begin{gather*}
 \deg u=\deg v=F,\quad \deg z=S,\quad
 \deg e_i=E_i,\quad \deg t_i=F-E_i,\\
 \deg\xi=2H-\sum_iE_i,\qquad
 \deg\eta=3H-\sum_iE_i.
\end{gather*}
For $s=1$, elimination of one base variable yields the six-generator hypersurface of Theorem~\ref{thm:rho3}. For $s\geq2$, elimination of both base variables yields a complete intersection with $2s+3$ minimal generators and $s-1$ minimal relations: $s-2$ relations in degree $F$, and one in degree $6H-2\sum_iE_i$. In particular, every surface satisfying these hypotheses is a Mori dream surface.
\end{theorem}
\begin{proof}
Write $R_0$ using \eqref{eq:multishift} and set $I_i=(\ell_i,X_1,y)$. The local tangent-cone argument of Theorem~\ref{thm:rho3} applies to each $I_i$, because the corresponding fiber is exactly $I_2$. Thus $\gr_{I_i}R_0$ is a domain with defining relation monic of degree two in the initial form of $y$, and $I_i$-adic order computes $\nu_i=\ord_{E_i}$.

Every element of $R_0$ has a unique normal form
\[
 h=\sum_{p\geq0,\ q=0,1}c_{pq}(u,v,z)X_1^py^q.
\]
The initial form of a nonzero such expression cannot vanish in $\gr_{I_i}R_0$: its degree in the initial form of $y$ is at most one, whereas that ring's defining relation is monic of degree two in that variable. Consequently there is the exact formula
\begin{equation}\label{eq:multiorder}
 \nu_i(h)=\min_{p,q}\bigl(\ord_{\ell_i}(c_{pq})+p+q\bigr).
\end{equation}
For a section in degree $aS+bF+\sum_i c_iE_i$, regularity along every exceptional curve therefore requires, term by term,
\begin{equation}\label{eq:divisibility}
 \prod_i\ell_i^{k_i}\mid c_{pq},\qquad
 k_i=\max\{0,-c_i-p-q\}.
\end{equation}
As in Proposition~\ref{prop:multifiltration} below, sections away from these curves extend to the normal Weierstrass model, so these conditions are sufficient as well as necessary.

Now substitute $X_1=\mathcal E\xi$, $y=\mathcal E\eta$, and $\ell_i=e_it_i$. A term satisfying \eqref{eq:divisibility}, after accounting for its exceptional degree, becomes a polynomial in the claimed generators: its exponent of $e_i$ is $c_i+p+q+k_i\geq0$, and its exponent of $t_i$ is $k_i$. This proves generation in every Picard degree. Substitution in \eqref{eq:multishift} and division by $\mathcal E^2$ give $G_s$.

We verify the kernel. The ring
\[
 B_0=\C[u,v,e_1,t_1,\ldots,e_s,t_s]/(e_it_i-\ell_i)
\]
is a domain: successively adjoining $e,t$ with $et=f\ne0$ to a domain embeds the resulting ring into the Laurent ring by $t=f/e$. Each $\ell_i$ remains nonzero during this construction. Since $G_s$ is monic in $\eta$, the final quotient is free with basis $1,\eta$ over $B_0[z,\xi]$. In particular every $e_i$ is a non-zero-divisor. Localizing at $\mathcal E$ recovers $R_0[e_1^{\pm1},\ldots,e_s^{\pm1}]$, so no further kernel can occur. The displayed equations form a regular sequence, by the same successive-domain argument followed by the monic equation.

For $s\geq2$, choose $\ell_1,\ell_2$ as base coordinates and eliminate $u,v$. The remaining $s-2$ base relations are quadratic relations among the products $e_it_i$. The last relation, as well as these quadrics, has no ordinary linear term. The homogeneous augmentation ideal modulo its square consequently has dimension $2s+3$, proving minimality of the generators. The regular sequence has length $s-1$, so its relations are minimal as well.
\end{proof}

\begin{corollary}\label{cor:multih0}
For $a\geq0$ and arbitrary integers $b,c_1,\ldots,c_s$,  $h^0\left(aS+bF+\sum_{k=1}^s c_kE_k\right)$ is equal to 
\begin{equation}\label{eq:multih0}
\begin{aligned}
\sum_{\substack{i\geq0,\ j\in\{0,1\}\\2i+3j\leq a}}
\max\bigg\{ b-\chi(2i+3j)-\sum_{k=1}^s \max\{0,-c_k-i-j\} +1, 0\bigg\}.
\end{aligned}
\end{equation}
For $a<0$ the dimension is zero.
\end{corollary}
\begin{proof}
In a fixed homogeneous normal-form term the exponent of $z$ is $a-2i-3j$. Its binary-form coefficient has degree $b-\chi(2i+3j)$. Condition \eqref{eq:divisibility} prescribes a factor of degree $\sum_k\max\{0,-c_k-i-j\}$. Counting the remaining binary forms proves the formula.
\end{proof}

\begin{corollary}\label{cor:endpoint}
If $s=2\chi+1$, the divisor of $\xi$ is an irreducible horizontal bisection $C$ with
\[
 [C]=2H-\sum_iE_i,\qquad C^2=-2,\qquad p_a(C)=\chi-2.
\]
In particular, the visible fiber components and the zero section no longer generate the effective cone, even though the Cox ring is finitely generated.
\end{corollary}
\begin{proof}
The exact order of $X_1$ along each $E_i$ is one, so division by $\mathcal E$ removes these components completely. The section $\xi$ does not vanish identically on any other fiber component: $x-rz^2$ is not identically zero on the corresponding Weierstrass cubic. It does not vanish along $S$ either. Its zero divisor therefore has only horizontal components and total fiber degree two. If reducible or nonreduced, it would have a component of fiber degree one, which is a section; the only section is $S$, already excluded. Thus it is an irreducible bisection. Intersection and adjunction give
\[
 C^2=4\chi-2s=-2,\qquad
 2p_a(C)-2=C^2+K_X\cdot C=-2+2(\chi-2).
\]
Every class in the cone generated by $S,E_i,F-E_i$ has nonnegative intersection with this distinct irreducible curve, whereas $C^2<0$. Hence its class is outside that cone.
\end{proof}

\begin{corollary}\label{cor:multicones}
Put $D_\xi=2H-\sum_iE_i$. Under the hypotheses of Theorem~\ref{thm:multiI2},
\[
 \Eff(X)=\operatorname{cone}(S,E_i,F-E_i,D_\xi).
\]
For $s\leq2\chi$, the last class is redundant. For $s=2\chi+1$, a class $D=aS+bF-\sum_i d_iE_i$ is nef precisely when
\[
 0\leq d_i\leq a/2\quad\text{for every }i,\qquad
 b\geq\chi a,\qquad b\geq\sum_i d_i.
\]
\end{corollary}
\begin{proof}
All generator degrees lie in the displayed cone: $F=E_i+(F-E_i)$, $H=S+\chi F$, and $\deg\eta=D_\xi+H$. Each displayed ray class is effective, so equality follows from the presentation. If $s\leq2\chi$, then
\[
 D_\xi=2S+\sum_i(F-E_i)+(2\chi-s)F.
\]
At the endpoint, duality gives the stated inequalities, using $D\cdot D_\xi=2b-2\sum_i d_i$ together with intersections with $S,E_i,F-E_i$.
\end{proof}
\subsection{Multigraded presentations of Cox rings}\label{subsec:multigraded}
The preceding presentations admit a common description in terms of exceptional valuations. 
We retain trivial Mordell--Weil group, so that the zero section, a fiber and the nonidentity fiber components form an integral Picard basis. 
Let $E_1,\ldots,E_r$ be those components and $f:X\to W$ the Weierstrass contraction. 
Choose a smooth fiber as the representative of $F$, so that $S,F$ are pullbacks of Cartier divisors whose supports avoid the singular points of $W$.
\begin{proposition}\label{prop:multifiltration}
For a homogeneous section $s\in (R_0)_{aS+bF}$, let $\nu_j(s)$ be its vanishing order along $E_j$ after pullback. Then, as a Picard-graded algebra,
\begin{equation}\label{eq:multifiltration}
 \Cox(X)\cong
 \bigoplus_{a,b,\,\boldsymbol c}
 \left\{s\in(R_0)_{aS+bF}:\nu_j(s)+c_j\geq0\text{ for every }j\right\}
 \boldsymbol t^{\boldsymbol c},
\end{equation}
where $\boldsymbol c\in\Z^r$, $\deg t_j=E_j$, and the degree of a summand is $aS+bF+\sum_jc_jE_j$. The algebra on the right is viewed inside $R_0[t_1^{\pm1},\ldots,t_r^{\pm1}]$.
\end{proposition}
\begin{proof}
A section of $aS+bF+\sum_jc_jE_j$ restricts away from the exceptional curves to a section of the corresponding Cartier divisor on the smooth open subset of $W$. 
The complement consists of finitely many points. 
Normality of $W$ extends this section across that complement. 
Its pullback defines a section of the desired divisor exactly when all exceptional orders are nonnegative after adding $c_j$. 
The condition is compatible with products, yielding the graded algebra identification.
\end{proof}
This formula is an exact description, but not in general a finite presentation. Theorem \ref{thm-I_2} computes it for several $I_2$ fibers under the stated hypotheses. 
For $I_n$ with $n\geq3$, there are $n-1$ exceptional valuations over a surface singularity of type $A_{n-1}$. 
One must determine their simultaneous vanishing ideals; ordinary powers of a single ideal no longer encode all of them. 
Saturated Rees constructions are a standard computational framework for such questions; see \cite{HKL16}.
\begin{problem}\label{prob:next}
Compute a finite multigraded presentation of \eqref{eq:multifiltration} first for one fiber of type $I_3$, and then for mixed configurations and for $I_2$ configurations beyond the interpolation range of Theorem~\ref{thm-I_2}, under the stated integral-Picard-basis hypotheses and appropriate finite-generation assumptions. In each case, establish that the proposed generators realize all simultaneous vanishing conditions, determine the kernel, and prove minimality.
\end{problem}
\begin{proof}[Proof of Theorem \ref{thm-I_2}]
It follows from Theorem \ref{thm:multiI2}, Corollaries 
\ref{cor:endpoint} and \ref{cor:multicones}.
\end{proof}
\section{Ample degrees and the cohomological criterion}
\label{sec:ample}
We prove Theorem~\ref{thm:no-ample} by a uniform cohomological argument. The analogue for extremal rational elliptic surfaces was considered in \cite[Proposition 2.1]{AGL16}. The additional step here is restriction to the zero section after applying Kawamata--Viehweg with the correct canonical twist.

\begin{proposition}\label{prop:two-fibers}
Let $D$ be an effective class with $D\not\sim F$. If $H^1(X,\OO_X(D-2F))=0$, then $\Cox(X)$ has no generator in degree $D$. In particular, the conclusion holds whenever $D-\chi F$ is nef and big.
\end{proposition}
\begin{proof}
Choose two distinct smooth fibers and use Proposition~\ref{prop:koszul}. Its multiplication map becomes
\[
 H^0(X,\OO_X(D-F))\otimes H^0(X,\OO_X(F))
 \longrightarrow H^0(X,\OO_X(D)).
\]
The stated cohomology vanishing implies surjectivity. For the second assertion,
\[
 (D-2F)-K_X=D-\chi F,
\]
so Kawamata--Viehweg vanishing applies under exactly the stated positivity hypothesis.
\end{proof}

\begin{lemma}\label{lem:nef-big-section}
Assume that $\pi$ has at least one reducible fiber and that $\NE(X)=\mathcal C_\pi$. If $A$ is a nef divisor satisfying $A\cdot S>0$ and $A\cdot F>0$, then $A$ is big.
\end{lemma}
\begin{proof}
On a smooth projective surface a nef class is pseudoeffective, and curve classes and divisor classes are identified numerically by the intersection pairing. Since $\mathcal C_\pi$ is closed, there is a numerical expression
\[
 A\equiv aS+\sum_{T\in\mathcal T}b_TT,
 \qquad a,b_T\geq0.
\]
Nefness gives $A^2\geq0$. If $A^2=0$, intersecting this expression with $A$ shows that $a=0$, because every term is nonnegative and $A\cdot S>0$. But then $A\cdot F=0$, since every $T$ is vertical. This is a contradiction. Thus $A^2>0$, which for a nef divisor on a surface is equivalent to bigness.
\end{proof}

\begin{theorem}\label{thm:uniform-vanishing}
Assume that $\pi$ has at least one reducible fiber and that $\NE(X)=\mathcal C_\pi$. For every ample integral divisor $D$, the divisor
\[
 A=D-\chi F-S
\]
is nef and big. Moreover,
\[
 H^1(X,\OO_X(D-2F-S))=H^1(X,\OO_X(D-2F))=0.
\]
In particular, multiplication by $H^0(X,\OO_X(F))$ is surjective in degree $D$.
\end{theorem}
\begin{proof}
The intersection identities $S^2=-\chi$, $S\cdot F=1$ and $F\cdot T=0$ give
\[
 A\cdot S=D\cdot S>0,\qquad
 A\cdot T=D\cdot T-S\cdot T\geq0
 \quad(T\in\mathcal T).
\]
Indeed, $D\cdot T\geq1$ by integrality and ampleness, while a section meets exactly one component of each reduced fiber, transversely, so $S\cdot T$ is zero or one. The cone hypothesis proves that $A$ is nef. If $F_i$ is a reducible fiber of type $I_{n_i}$, then
\[
 D\cdot F=\sum_{T\subset F_i}D\cdot T\geq n_i\geq2.
\]
Consequently $A\cdot F=D\cdot F-1>0$, and Lemma~\ref{lem:nef-big-section} proves that $A$ is big.

Since $K_X+A\sim D-2F-S$, Kawamata--Viehweg vanishing gives
\[
 H^1(X,\OO_X(D-2F-S))=0.
\]
The restriction sequence
\[
 0\longrightarrow\OO_X(D-2F-S)
 \longrightarrow\OO_X(D-2F)
 \longrightarrow\OO_S(D-2F)
 \longrightarrow0
\]
has last term of degree $D\cdot S-2\geq-1$ on $S\cong\PP^1$. Its first cohomology therefore vanishes as well. The long exact sequence proves $H^1(X,\OO_X(D-2F))=0$. Applying the Koszul sequence of two distinct smooth fibers, as in Proposition~\ref{prop:two-fibers}, gives the asserted surjectivity.
\end{proof}

\begin{proof}[Proof of Theorem~\ref{thm:no-ample}]
Lemma~\ref{lem:MW} gives trivial Mordell--Weil group and a free Picard group. If there is a reducible fiber, the cone theorem of \cite[Theorem 1.1]{LLL26} gives $\NE(X)=\mathcal C_\pi$, and Theorem~\ref{thm:uniform-vanishing} proves part~(i).

If there are no reducible fibers, Shioda--Tate gives $\rho(X)=2$. By Corollary~\ref{cor:rho2cones}, an ample integral class has the form $D=aS+bF$ with $a\geq1$ and $b\geq a\chi+1$. Corollary~\ref{cor:section-fiber-cohomology} gives
\[
 H^1(X,\OO_X(aS+(b-2)F))=0,
\]
because $b-2\geq a\chi-1$. The same fiber-pencil argument proves part~(i).

For part~(ii), put $d=D\cdot S\geq1$. If there is a reducible fiber, then for $0\leq k\leq d$,
\[
 (D-kF)\cdot S=d-k,\qquad
 (D-kF)\cdot T=D\cdot T>0\quad(T\in\mathcal T).
\]
Thus $D-kF$ is ample for $k<d$ and nef for $k=d$, by the description of the closed cone of curves and Kleiman's criterion. The same conclusion in Picard number two follows by testing against $S$ and $F$. Applying part~(i) successively to $D,D-F,\ldots,D-(d-1)F$ yields a surjective multiplication map with $d$ factors $H^0(X,\OO_X(F))$. Commutativity makes this map factor through their symmetric power, giving part~(ii).

Finally, an ample class is not $F$. Every nonzero term in the image of the multiplication map in part~(i) is therefore a product of sections of two nonzero effective degrees. In the notation of Section~\ref{sec-nef-degree}, $(R_+/R_+^2)_D=0$. If $R$ is finitely generated, Proposition~\ref{prop:negative-nef} shows that every minimal generator other than a defining section of a negative curve has nef degree; since the interior of the nef cone is the ample cone, such a degree lies on its boundary. This proves part~(iii).
\end{proof}
\begin{corollary}\label{cor:explicit-no-ample}
Assume either that $\rho(X)=2$, or that all $s\geq1$ reducible fibers are of type $I_2$, $\MW(\pi)$ is finite, and $s\leq2\chi+1$. Then no ample divisor class is the degree of a minimal homogeneous generator of $\Cox(X)$.
\end{corollary}
\begin{proof}
The Picard-number-two case follows from Theorem~\ref{thm:no-ample}, or directly from the degrees $S,F,2H,3H$ in Theorem~\ref{thm:rho2}.
For $s\geq1$, Theorem~\ref{thm:multiI2} gives the minimal degrees
\[
 S,\quad E_i,\quad F-E_i,\quad
 2H-\sum_iE_i,\quad 3H-\sum_iE_i,
\]
with one additional degree $F$ when $s=1$. The first three types are negative-curve classes, and $F$ is not ample. Each of the last two degrees has intersection zero with $S$, and hence is not ample either. The degrees of the indecomposables $R_+/R_+^2$ are intrinsic, so this holds for every minimal homogeneous generating set. In particular, the argument includes $s=2\chi+1$, which is outside the scope of Theorem~\ref{thm:no-ample}.
\end{proof}

\begin{remark}\label{rem:canonical-twist}
More generally, for disjoint curves $B_1,B_2$, the Kawamata--Viehweg hypothesis relevant to $H^1(D-B_1-B_2)$ is that
\[
 D-B_1-B_2-K_X
   =D-B_1-B_2-(\chi-2)F
\]
be nef and big. Nefness and bigness of $D-B_1-B_2$ alone do not justify this vanishing. Likewise, ampleness of $D$ does not imply that $D-\chi F$ is nef. Theorem~\ref{thm:uniform-vanishing} instead applies vanishing to $D-2F-S$ and then uses restriction to $S$. The group needed for the fiber pencil is $H^1(D-2F)$; no vanishing of $H^1(D-\chi F)$ is asserted or required.
\end{remark}

\begin{example}
\label{ex:nef-big-h1}
Suppose that all reducible fibers are $I_2$, their number is $s=2\chi$, and $\MW(\pi)$ is finite. Set
\[
 L=4H+F-2\sum_{i=1}^{2\chi}E_i.
\]
By Corollary~\ref{cor:multicones}, nefness can be checked against $S,E_i,F-E_i$. We have
\[
 L\cdot S=1,\qquad L\cdot E_i=4,\qquad
 L\cdot(F-E_i)=0,\qquad L^2=8.
\]
Thus $L$ is nef and big. The four possible normal-form terms in Corollary~\ref{cor:multih0} have pole orders $0,2,3,4$; they contribute respectively $2,2,\pos{2-\chi},2$. Hence $h^0(L)=6$. Since $L\cdot F=4$, Serre duality gives $h^2(L)=0$, while Riemann--Roch gives
\[
 \chi(\OO_X(L))=\chi+\tfrac12\bigl(8-4(\chi-2)\bigr)=8-\chi.
\]
Consequently
\[
 h^1(X,\OO_X(L))=\chi-2>0.
\]
The $8I_2+32I_1$ example with $\chi=4$ in \cite[Theorem 3.8]{LL26} realizes these hypotheses, and gives $h^1(L)=2$. Also $(L-K_X)\cdot S=3-\chi$; for that example this intersection is negative, so the relevant vanishing hypothesis fails.
\end{example}

\begin{remark}
The Mori dream property alone must not be used to identify the Mori cone with the cone generated by $S$ and the reducible fiber components. Under the finite Mordell--Weil hypotheses, the cited cone theorem requires $\delta(\pi)\leq\chi$. At the endpoint $s=2\chi+1$, Corollaries~\ref{cor:endpoint} and~\ref{cor:multicones} exhibit the additional negative bisection, although Theorem~\ref{thm:multiI2} proves finite generation.
\end{remark}

\begin{corollary}\label{cor:direct-fiber-criterion}
Assume that $\MW(\pi)$ is finite and $\delta(\pi)\leq\chi$. For an ample integral divisor $D$, the divisor $D-\chi F$ is nef and big if and only if $D\cdot S\geq\chi$.
\end{corollary}
\begin{proof}
Put $B=D-\chi F$. Then $B\cdot S=D\cdot S-\chi$, while $B$ has strictly positive intersection with every vertical generator of the Mori cone. Here those generators are the reducible fiber components if any exist, and $F$ in Picard number two. The descriptions of the cones in \cite[Theorem 1.1]{LLL26} and Corollary~\ref{cor:rho2cones} therefore show that $B$ is nef exactly when $D\cdot S\geq\chi$.

In that case $B$ is pseudoeffective. If $B^2=0$, expressing $B$ as a nonnegative combination of the cone generators and intersecting with $B$ forces every vertical coefficient to vanish. Thus $B\equiv aS$. Since $a=B\cdot F=D\cdot F>0$, this gives $B^2=-\chi a^2<0$, a contradiction. Hence $B$ is also big.
\end{proof}

\begin{remark}
When $1\leq D\cdot S<\chi$, the divisor $D-\chi F$ is not nef, but Theorem~\ref{thm:no-ample} still gives $H^1(X,\OO_X(D-2F))=0$ and excludes $D$ as a generator degree. Thus small intersection with the zero section produces no exceptional ample degrees in the cone range.
\end{remark}

\begin{example}\label{ex:endpoint-vanishing}
At the interpolation endpoint, the explicit cohomology formula gives a different behavior. Suppose that $X$ satisfies the hypotheses of Theorem~\ref{thm:multiI2} with $s=2\chi+1$, and put
\[
 D=5H+F-2\sum_{i=1}^{2\chi+1}E_i.
\]
Its intersections with $S,E_i,F-E_i$, and the additional bisection $C$ of Corollary~\ref{cor:endpoint} are, respectively,
\[
 1,\qquad4,\qquad1,\qquad2\chi-2.
\]
They are all positive, so Corollary~\ref{cor:multicones} and Kleiman's criterion show that $D$ is ample. For $L=D-2F$, formula~\eqref{eq:multih0} gives
\[
 h^0(L)=(\chi-2)+(\chi-1)+\chi=3\chi-3.
\]
Moreover $L^2=9\chi-18$, $K_X\cdot L=5(\chi-2)$ and $h^2(L)=0$, the latter because $(K_X-L)\cdot F=-5$. Riemann--Roch yields
\[
 \chi(\OO_X(L))=3\chi-4,\qquad h^1(L)=1.
\]
Thus, for this fiber configuration, uniform vanishing fails although $D$ is not a generator degree by Corollary~\ref{cor:explicit-no-ample}.
\end{example}
\section{A finite criterion for Cox generators}
\label{sec:semistable-case}
\subsection{The finite criterion}
\begin{proof}[Proof of Theorem \ref{thm:criterion}]
We prove finiteness of the proposed set, generation, and minimality,
in that order.

\medskip\noindent\textbf{ Step 1. Finiteness of the monoid and the candidate degrees.}
Since $\MW(\pi)=0$, the classes $S,F,T_{i,j}$ with $j\ge1$ form an
integral basis of $\Pic(X)$. By Theorem~\ref{cone-thm}, $M$ is the monoid of lattice points in a pointed rational polyhedral cone. 
It is therefore finitely generated by Gordan's lemma.

We first check that $0$ is its only element of fiber degree zero, and
that $H$ is its only element of fiber degree one. For $A\in M$, set
$m=A\cdot F$ and write
\[
A=mH-\sum_{i=1}^s\sum_{j=1}^{n_i-1}b_{i,j}T_{i,j},
\qquad b_{i,0}=b_{i,n_i}=0.
\]
The integers $w_{i,j}=b_{i,j}-b_{i,j-1}$ satisfy
\[
w_{i,1}\ge\cdots\ge w_{i,n_i},\qquad
w_{i,1}-w_{i,n_i}\le m,\qquad \sum_jw_{i,j}=0,
\]
because $A$ has nonnegative intersection with each fiber component.
If $m=0$ or $m=1$, all these integers are zero: a list of integers of
width at most one and mean zero must be identically zero.
Thus $A=mH$. Moreover, $H$ is nef by Theorem~\ref{cone-thm}.

Consequently, the Hilbert basis of $M$ has the form
$\{H\}\cup\mathcal Q$, with $Q\cdot F\ge2$ for $Q\in\mathcal Q$.
The monoid $M^*=M\setminus\{H\}$ is generated by the finite list
\begin{equation}\label{eq:finite-list}
\mathcal Q,\qquad 2H,\qquad 3H,\qquad H+Q\quad(Q\in\mathcal Q).
\end{equation}
Indeed, write an element of $M$ as $aH+\sum Q_k$.
Every integer $a\ge2$ is a nonnegative combination of $2$ and $3$;
if $a=1$ and some $Q_k$ occurs, group $H$ with that $Q_k$.
This proves~\eqref{eq:finite-list}.
Adding the nonidentity fiber components shows that $P$ is finitely
generated. It is pointed: intersection with an ample integral divisor
is strictly positive on every nonzero summand in its definition.
Its indecomposable elements consequently form a finite set and generate
$P$.

If $G\in\mathcal G$, an expression
$G=A+\sum a_{i,j}T_{i,j}$ with $A\in M^*$ must have $A\ne0$.
Indecomposability forces all $a_{i,j}=0$, so $G=A\in M^*$.
This proves the first assertion of the theorem. It also shows that every
element of $P$ is a sum of elements of $\mathcal G$ and nonidentity
fiber components: the elements of fiber degree zero in $P$ are precisely
nonnegative integral combinations of those components.

\medskip\noindent\textbf{Step 2. Sections of vertical degree.}
Let $R$ be the subalgebra of $\Cox(X)$ generated by the displayed finite
set. Distinct reducible fibers give nonproportional sections
$p_i\in H^0(X,\mathcal O_X(F))$. Thus any two of them form a basis,
and the additional sections in the statement ensure
$H^0(X,\mathcal O_X(F))\subset R$ for every $s\ge1$.
Since $\pi_*\mathcal O_X=\mathcal O_{\mathbb P^1}$,
\[
H^0(X,\mathcal O_X(kF))
=H^0(\mathbb P^1,\mathcal O_{\mathbb P^1}(k))\subset R
\qquad(k\ge0).
\]
Here the inclusion follows because the homogeneous coordinate ring of
$\mathbb P^1$ is generated in degree one.

If $D$ is effective and $D\cdot F=0$, the zero divisor of any nonzero
section of $D$ is vertical. It is a sum of components of reducible fibers
and whole irreducible fibers. The defining section of each whole fiber
belongs to $H^0(X,\mathcal O_X(F))$. Therefore that section of $D$ is,
up to a nonzero scalar, a product of elements of $R$.
This proves $H^0(X,\mathcal O_X(D))\subset R$ in every vertical degree.

\medskip\noindent\textbf{Step 3. Generation in all remaining degrees.}
Fix an ample integral divisor $A$ and induct on the positive integer
$A\cdot D$ for nonzero effective classes $D$.
If $D$ is not nef, Theorem~\ref{cone-thm} gives a curve
$C\in\{S,T_{i,j}\}$ such that $D\cdot C<0$; indeed $D\cdot F\ge0$.
Every effective representative of $D$ contains $C$, so
\[
H^0(X,\mathcal O_X(D))=s_C H^0(X,\mathcal O_X(D-C)).
\]
The induction hypothesis applies to $D-C$ whenever its space of sections is nonzero.

We may therefore assume that $D$ is nef. 
 Fiber degree zero was treated
above. 
If $D\cdot F=1$, every effective representative of $D$ has exactly one horizontal component, with multiplicity one and degree one over the base. 
Its generic point defines a rational point of the generic elliptic
fiber. 
Since $\MW(\pi)=0$, this component is $S$. 
Thus $S$ is fixed, and the same induction applies with $C=S$.

Suppose now that $D\cdot F\ge2$. Put
\[
b=D\cdot S\ge0,\qquad D_0=D-bF.
\]
Then $D_0\cdot S=0$, $D_0\cdot F=D\cdot F$, and
$D_0\cdot T_{i,j}=D\cdot T_{i,j}\ge0$.
By Theorem~\ref{cone-thm}, $D_0$ is nef, hence $D_0\in M^*\subset P$.
Decompose $D_0$ into elements of $\mathcal G$ and nonidentity components.
The product of the corresponding sections is an element
\[
g\in R\cap H^0(X,\mathcal O_X(D_0))
\quad\text{with}\quad g|_S\ne0,
\]
because every nonidentity component is disjoint from $S$.
Since $D_0\cdot S=0$, its restriction $g|_S$ is a nonzero constant
after trivializing $\mathcal O_S(D_0)$.

Since $\pi|_S\colon S\to\mathbb P^1$ is an isomorphism, restriction gives
\[
H^0(X,\mathcal O_X(bF))\simeq H^0(S,\mathcal O_S(bF)).
\]
Multiplication by $g|_S$ identifies the space on the right with
$H^0(S,\mathcal O_S(D))$.
For any section of $D$, we can therefore subtract an element of
$gH^0(X,\mathcal O_X(bF))\subset R$ having the same restriction to $S$.
The difference vanishes on $S$, and the restriction sequence gives
\begin{equation}\label{eq:reduction}
\begin{split}
H^0(X,\mathcal O_X(D))
={}&gH^0(X,\mathcal O_X(bF))\\
&+s_SH^0(X,\mathcal O_X(D-S)).
\end{split}
\end{equation}
The second summand belongs to $R$ by induction, since
$A\cdot(D-S)<A\cdot D$. The first already belongs to $R$, including
when $b=0$. Thus $R=\Cox(X)$.
This proves finite generation without any prior semiampleness assumption.

\medskip\noindent\textbf{Step 4. Minimality of the horizontal generators.}
We need one observation. If $B\cdot S=0$ and
$t\in H^0(X,\mathcal O_X(B))$ has $t|_S\ne0$, then $B\in P$.
Indeed, the zero divisor of $t$ is disjoint from $S$.
If $B$ is not nef, choose a curve with negative intersection as in
Step~3 and divide out its defining section. That curve is a fixed
component of the zero divisor of $t$; it cannot be $S$ or any $T_{i,0}$,
since those meet $S$. It is therefore a nonidentity component.
Repeating this process terminates by ample degree and gives
\[
B=B'+\sum_{i,j\ge1}a_{i,j}T_{i,j},\qquad B'\in M,
\]
together with a section of $B'$ nonzero on $S$.
If $B'\cdot F=0$, then $B'=0$ by Step~1.
If $B'\cdot F=1$, then $B'=H$, but every section of a class of fiber
degree one vanishes on $S$, by Step~3. This is impossible.
Hence $B'\in M^*$, proving the observation.

Fix $G\in\mathcal G$. If $g_G$ were a sum of products of at least two
sections of nonzero effective degrees, at least one such product would restrict
nontrivially to $S$. Each of its factors would then restrict nontrivially
to $S$ and have nonnegative intersection with $S$.
Their degrees sum to $G$, which is orthogonal to $S$, so every factor
degree is orthogonal to $S$. The observation places all these nonzero
degrees in $P$. A product of at least two factors would therefore
decompose $G$ in $P$, contradicting its definition.
Thus every $g_G$ is indispensable in the proposed generating set.

\medskip\noindent\textbf{Step 5. Minimality of the remaining generators.}
Each of $S$ and $T_{i,j}$ is an irreducible curve of negative square.
Its complete linear system consists only of that curve. A product of
sections of two nonzero effective degrees cannot have this irreducible
zero divisor, so its defining section is indispensable.

Finally, a monomial of degree $F$ cannot involve a horizontal generator,
since such a generator has positive intersection with $F$.
A monomial made only from fiber-component sections has degree $F$
precisely when its zero divisor is one of the reducible fibers; these
monomials are the $p_i$, up to scalars.
If $s=1$, one section independent of $p_1$ is necessary.
If $s\ge2$, no additional section is needed.
This proves minimality of the entire set.

The surface is smooth and projective, with $q(X)=0$ and free finitely
generated Picard group. The standard Cox-ring characterization of Mori
dream spaces therefore gives the last assertion.
\end{proof}
\subsection{Nonvanishing in the numerical range}
The following theorem supplies every nonvanishing section required by Theorem~\ref{thm:criterion}. 
\begin{theorem}
\label{thm-2chi+4}
Assume that $s\ge1$, $\operatorname{MW}(\pi)=0$, and $N(\pi)\le 2\chi+4$.
Then, for every nef integral divisor $G$ satisfying
$G\cdot S=0$ and $G\cdot F\ge2$, the restriction map
\[
H^0(X,\mathcal O_X(G))
\longrightarrow H^0(S,\mathcal O_S(G))\simeq\mathbb C
\]
is surjective.
Consequently, no nef and big integral divisor $D$ with $D\cdot S>0$
is a necessary degree of $\Cox(X)$.
\end{theorem}
\begin{proof}
From the short exact sequence
\[
0\longrightarrow\mathcal O_X(G-S)
\longrightarrow\mathcal O_X(G)
\longrightarrow\mathcal O_S(G)
\longrightarrow0,
\]
the restriction map on global sections is surjective whenever
$H^1(X,\mathcal O_X(G-S))=0$. We prove this vanishing except in one
case, where the required section is a power of a section of smaller degree.
We then deduce the assertion about $\operatorname{Cox}(X)$.

\medskip\noindent\textbf{Step 1: Coordinates in $\operatorname{Pic}(X)$ and nefness.}
Put $H=S+\chi F$ and $m=G\cdot F\ge 2$. Since $S^2=-\chi$, we have
$H\cdot S=0$ and $H\cdot F=1$. Because $\operatorname{MW}(\pi)=0$,
the classes $S$, $F$, and the non-identity components $T_{i,j}$
($1\le j\le n_i-1$) of the reducible fibers form an integral basis of
$\operatorname{Pic}(X)$. Denote by $T_{i,0}$ the identity component of the
$i$-th reducible fiber (the unique component meeting $S$), and order
$T_{i,0},T_{i,1},\ldots,T_{i,n_i-1}$ cyclically. Since $G\cdot S=0$ and
$G\cdot F=m$, we can write
\[
G=mH-\sum_{i=1}^s\sum_{j=1}^{n_i-1}b_{i,j}T_{i,j}
\qquad (b_{i,j}\in\mathbb Z).
\]
Set $b_{i,0}=b_{i,n_i}=0$ and define the first-difference vector
$w_i=(w_{i,1},\ldots,w_{i,n_i})\in\mathbb Z^{n_i}$ by
\[
w_{i,j}=b_{i,j}-b_{i,j-1}\qquad (1\le j\le n_i).
\]
Using the intersection matrix of the $I_{n_i}$ fiber, we obtain
(the two components of an $I_2$ fiber intersect twice)
\[
G\cdot T_{i,j}=w_{i,j}-w_{i,j+1}\quad (1\le j\le n_i-1),
\qquad
G\cdot T_{i,0}=m-(w_{i,1}-w_{i,n_i}).
\]
Thus $G$ is nef on the fibers if and only if, for each $1\le i\le s$,
the vector $w_i$ satisfies:
\begin{equation}\label{eq:nef}
w_{i,1}\ge w_{i,2}\ge\cdots\ge w_{i,n_i},\qquad
w_{i,1}-w_{i,n_i}\le m,\qquad
\sum_{j=1}^{n_i} w_{i,j}=0.
\end{equation}
We call $w_{i,1}-w_{i,n_i}$ the \emph{width} of the sorted vector $w_i$.

\medskip\noindent\textbf{Step 2: Cohomological reduction rules.}
Let $L$ be an integral divisor on $X$.
\begin{enumerate}
\item[(a)] If $C\simeq\mathbb P^1$ is a fiber component with $L\cdot  C=-1$,
then $$H^0(C,\mathcal O_C(L))=H^1(C,\mathcal O_C(L))=0,$$
 so restriction to $C$ yields an isomorphism
\[
H^1(X,\mathcal O_X(L-C))\xrightarrow{\;\sim\;} H^1(X,\mathcal O_X(L)).
\]
\item[(b)] If $L\cdot S\ge -1$, then $H^1(S,\mathcal O_S(L))=0$, so
restriction to $S$ yields a surjection
\[
H^1(X,\mathcal O_X(L-S))\twoheadrightarrow H^1(X,\mathcal O_X(L)).
\]
\end{enumerate}

\medskip\noindent\textbf{Step 3: The width-reduction algorithm.}
Starting with $L_0=G$, we construct a sequence of divisors
$L_1,\ldots,L_{m-1}$ with $L_k\cdot F=m-k$, each satisfying the nefness
conditions \eqref{eq:nef} with $m$ replaced by $m-k$, as follows.

To pass from $L_{k-1}$ (fiber degree $a+1$, where $a=m-k$) to $L_k$, we first
subtract $S$. Since $S\cdot T_{i,0}=1$ and $S\cdot T_{i,j}=0$ for $j\ge 1$,
subtracting $S$ leaves the vectors $w_i$ unchanged, while the intersection
with $T_{i,0}$ becomes $a-(w_{i,1}-w_{i,n_i})$.
If a vector $w_i$ has width $a+1$, then the current divisor has intersection
$-1$ with $T_{i,0}$. We then restore width at most $a$ by the following loop:
\begin{enumerate}
\item Subtract $T_{i,0}$: since $T_{i,0}\sim F-\sum_{j=1}^{n_i-1}T_{i,j}$,
this decreases $b_{i,j}$ by $1$ for all $1\le j\le n_i-1$,
which decreases $w_{i,1}$ by $1$ and increases $w_{i,n_i}$ by $1$.
\item Restore the non-increasing order of $w_i$: any adjacent inversion
created by step 1 has the form $(w_{i,j},w_{i,j+1})=(c,c+1)$, which means that
the current divisor has intersection $w_{i,j}-w_{i,j+1}=-1$ with $T_{i,j}$.
Subtracting $T_{i,j}$ swaps $w_{i,j}$ and $w_{i,j+1}$. Repeating these swaps
sorts the vector $w_i$; the net effect of steps 1 and 2 on the sorted vector
$w_i$ is to decrease a maximal entry by $1$ and increase a minimal entry by $1$.
\item Repeat steps 1 and 2 until $w_i$ has width at most $a$.
\end{enumerate}
Each maximum--minimum operation decreases $\sum_j w_{i,j}^2$ by
$2a>0$, while the swaps leave it unchanged. Thus the loop terminates. Because every fiber component $T_{i,j}$ ($0\le j\le n_i-1$) is
subtracted when the current divisor has intersection $-1$ with it, Rule (a)
gives isomorphisms
\[
\begin{aligned}
H^1(X,\mathcal O_X(G-S))&\simeq H^1(X,\mathcal O_X(L_1)),\\
H^1(X,\mathcal O_X(L_k-S))&\simeq H^1(X,\mathcal O_X(L_{k+1})),
\end{aligned}
\]
where the second isomorphism holds for $1\le k\le m-2$.
Let $r_k$ be the total number of copies of the identity components $T_{i,0}$
subtracted in passing from $G$ to $L_k$. Since $S^2=-\chi$ and $T_{i,0}\cdot S=1$,
we have
\[
L_k\cdot S = k\chi - r_k.
\]
By Rule (b), we obtain surjections
\[
H^1(X,\mathcal O_X(L_{k+1}))\twoheadrightarrow
H^1(X,\mathcal O_X(L_k))\qquad (1\le k\le m-2)
\]
provided that $L_k\cdot S\ge -1$ for all $1\le k\le m-2$, that is,
\begin{equation}\label{eq:counts}
r_k\le k\chi+1\qquad (1\le k\le m-2),
\end{equation}
It therefore suffices to verify \eqref{eq:counts} and
$H^1(X,\mathcal O_X(L_{m-1}))=0$.
When $m=2$, condition \eqref{eq:counts} is empty.

\medskip\noindent\textbf{Step 4: Combinatorial bound on $r_k$.}
Consider a nonzero sorted integer vector $w=(w_1,\ldots,w_n)$ with
$\sum_j w_j=0$, maximum $p>0$, minimum $-q<0$, and width $\ell=p+q$.
To reduce its width to $\ell-h\ge 1$ via maximum--minimum operations, choose
\[
u\in\left\{
\left\lfloor\frac{hq}{\ell}\right\rfloor,
\left\lceil\frac{hq}{\ell}\right\rceil
\right\}
\]
and consider the target interval $[-q+u,\,p-h+u]$ of length $\ell-h$.
Since $0\le u\le q$ and $0\le h-u\le p$, this interval contains $0$.
Define the total deficits outside this interval by
\[
A_u=\sum_{j=1}^n(-q+u-w_j)_+,\qquad
B_u=\sum_{j=1}^n(w_j-p+h-u)_+,
\]
where $x_+=\max\{x,0\}$. As long as the width of $w$ exceeds $\ell-h$, each
maximum--minimum operation decreases $\max\{A_u,B_u\}$ by $1$: if both $A_u$
and $B_u$ are positive, both decrease by $1$; if only one is positive, it
decreases by $1$ while the other remains $0$ (since the interval contains $0$
and $\sum_j w_j=0$, the entries cannot all lie on one side of $0$ unless
they are all $0$). Thus the number of operations is at most $\max\{A_u,B_u\}$.

Since the function $x\mapsto(-q+u-x)_+$ is convex on $[-q,p]$, it is bounded
above by the secant line connecting $(-q,u)$ and $(p,0)$, namely
$\frac{u}{\ell}(p-x)$. Summing over $j$ and using $\sum_j w_j=0$ gives
\[
A_u\le\frac{npu}{\ell},\qquad\text{and similarly}\qquad
B_u\le\frac{nq(h-u)}{\ell}.
\]
Setting $\alpha=p/\ell$, $\beta=q/\ell$ (so $\alpha+\beta=1$ and $\alpha\beta\le 1/4$),
and $f=h\beta-\lfloor h\beta\rfloor$, the minimum over the two choices of $u$ yields
\[
\min_u\max\{A_u,B_u\}
\le n\bigl(h\alpha\beta+\min\{f\beta,(1-f)\alpha\}\bigr)
\le n(h+1)\alpha\beta
\le\frac{n(h+1)}{4}.
\]
At stage $k$, each vector $w_i$ of initial width $\ell_i\le m$ is reduced to
width at most $m-k$, which requires a width reduction of
$h_i=\max\{0,\ell_i-m+k\}\le k$. If $h_i=0$, no operation is needed.
Otherwise, the cumulative maximum--minimum operations are the same as
those obtained by reducing the initial vector directly to width $m-k$:
the intermediate stages merely pause the same algorithm.
Summing over all reducible fibers gives
\[
r_k\le\left\lfloor\frac{N(k+1)}{4}\right\rfloor.
\]
Since $N\le 2\chi+4$ and $\chi\ge 3$, we have
$\lfloor N(k+1)/4\rfloor\le k\chi+1$ for all $k\ge 2$ (for $\chi=3$ and $k=2$,
note that $\lfloor 10\cdot 3/4\rfloor=7=2\chi+1$).

For $k=1$, a vector $w_i$ of width strictly less than $m$ requires no operation,
whereas a vector of width $m$ requires $\min\{\#\{j:w_{i,j}=p\},\,\#\{j:w_{i,j}=-q\}\}\le\lfloor n_i/2\rfloor$
operations. Hence
\[
r_1\le\sum_{i=1}^s\left\lfloor\frac{n_i}{2}\right\rfloor\le\frac{N}{2}\le\chi+2,
\]
so $r_1\le\chi+1$ holds except in the single extremal case where $N=2\chi+4$,
every $n_i$ is even, and every vector $w_i$ consists of $n_i/2$ copies of $m/2$
and $n_i/2$ copies of $-m/2$.

\medskip\noindent\textbf{Step 5: Vanishing at the final stage and the extremal case.}
At stage $k=m-1$, each reduced vector has integer entries, sum zero, and width
at most $m-(m-1)=1$, so it is identically zero. Since each maximum--minimum operation decreases
the sum of the positive entries $\sum_j(w_{i,j})_+=\max_j b_{i,j}$ by $1$, the
total number of operations performed up to stage $m-1$ is
$r_{m-1}=\sum_i\max_j b_{i,j}$. Therefore, the coefficients of the nonidentity components vanish in
$L_{m-1}$, and
\[
L_{m-1}=S+BF,\qquad B=m\chi-\sum_{i=1}^s\max_j b_{i,j}.
\]
For each nonzero initial vector, with maximum $p$ and minimum $-q$,
the secant-line estimate for $x\mapsto x_+$ gives
\[
\max_j b_{i,j}=\sum_{j=1}^{n_i}(w_{i,j})_+\le\frac{n_i pq}{p+q}\le\frac{mn_i}{4},
\]
For a zero vector the same inequality is immediate. Hence
\[
B\ge m\left(\chi-\frac{N}{4}\right)\ge m\left(\frac{\chi}{2}-1\right)>0.
\]
For completeness, on every fiber $C$ the line bundle
$\mathcal O_C(-S)$ has no section: on the component meeting $S$ its
degree is $-1$, and vanishing then propagates through the remaining
components, on which its degree is zero. Since $\omega_C\simeq\mathcal O_C$,
duality gives $H^1(C,\mathcal O_C(S))=0$ and hence
$h^0(C,\mathcal O_C(S))=1$. Cohomology and base change now yield
$R^1\pi_*\mathcal O_X(S)=0$, and the defining section of $S$ identifies
$\pi_*\mathcal O_X(S)$ with $\mathcal O_{\mathbb P^1}$.
The Leray spectral sequence therefore gives
\[
H^1(X,\mathcal O_X(L_{m-1}))=H^1(X,\mathcal O_X(S+BF))
\simeq H^1(\mathbb P^1,\mathcal O_{\mathbb P^1}(B))=0.
\]
Thus $H^1(X,\mathcal O_X(G-S))=0$ outside the extremal case.

In the extremal case, $m$ is necessarily even and $G=\frac{m}{2}Q$, where
\[
Q=2H-\sum_{i=1}^s\sum_{j=1}^{n_i-1}\min\{j,n_i-j\}T_{i,j}.
\]
Applying Step 3 to $Q$ (where $m=2$, so \eqref{eq:counts} is vacuous) yields
$L_1=S+(\chi-2)F$, so $H^1(X,\mathcal O_X(Q-S))\simeq H^1(\mathbb P^1,\mathcal O_{\mathbb P^1}(\chi-2))=0$.
Hence $Q$ admits a global section not vanishing on $S$, and its $(m/2)$-th
tensor power gives the desired section of $G$.

\medskip\noindent\textbf{Step 6: Application to $\operatorname{Cox}(X)$.}
Let $D$ be a nef and big integral divisor with $b=D\cdot S>0$.
Then $D\cdot F>0$. Since
\[
\delta(\pi)=\sum_{i=1}^s\frac{\lfloor n_i^2/4\rfloor}{n_i}\le\frac{N}{4}<\chi,
\]
the Mori cone $\overline{\operatorname{NE}}(X)$ is generated by $S$, $F$, and the
fiber components $T_{i,j}$. Consequently, $G:=D-bF$ is nef and satisfies
$G\cdot S=0$.
If $D\cdot F=1$, then $\operatorname{MW}(\pi)=0$ implies that $S$ is a fixed
component of $|D|$, so $H^0(X,\mathcal O_X(D))=s_S H^0(X,\mathcal O_X(D-S))$
(where $s_S\in H^0(X,\mathcal O_X(S))$ defines $S$), and $D$ is not a necessary
generator degree.
If $D\cdot F\ge 2$, then $G\cdot F\ge 2$, so by the first part of the theorem
there exists $g\in H^0(X,\mathcal O_X(G))$ whose restriction to $S$ is a nonzero
constant. Since $\mathcal O_S(D)\simeq\mathcal O_S(bF)\simeq\mathcal O_{\mathbb P^1}(b)$,
the restriction map $H^0(X,\mathcal O_X(bF))\to H^0(S,\mathcal O_S(bF))$ is
an isomorphism, which implies
\[
\begin{split}
H^0(X,\mathcal O_X(D))
={}&g\,H^0(X,\mathcal O_X(bF))\\
&+s_S\,H^0(X,\mathcal O_X(D-S)).
\end{split}
\]
Since both summands consist of products of sections of strictly smaller nonzero
effective degrees, $D$ is not a necessary degree of $\operatorname{Cox}(X)$.
\end{proof}
\begin{proof}[Proof of Theorem \ref{thm:bound}]
Since $\chi\ge3$, the numerical assumption gives
\[
\delta(\pi)\le\frac{N(\pi)}4\le\frac\chi2+1<\chi.
\]
Thus the cone criterion applies, and Lemma~\ref{lem:MW} gives $\MW(\pi)=0$. Every $G\in\mathcal G$ is nef,
orthogonal to $S$, and has $G\cdot F\ge2$.
Theorem~\ref{thm-2chi+4} supplies a section $g_G$ nonzero on $S$.
Theorem~\ref{thm:criterion} proves finite generation and minimality;
the last assertion is also part of Theorem~\ref{thm-2chi+4}.
\end{proof}
\subsection{One reducible fiber: explicit generator degrees}
\label{subsec:single-fiber}
The general criterion also recovers concrete generator degrees without a
separate cohomological argument for each fiber type.

\begin{corollary}\label{cor:one-fiber-application}
Assume that $\MW(\pi)$ is finite and that the unique reducible fiber is
$F_0=T_0+\cdots+T_{n-1}$ of type $I_n$, with cyclically ordered components,
$S\cdot T_0=1$, and $3\leq n\leq2\chi+4$.
Then $\Cox(X)$ has a minimal homogeneous generating set consisting of
the defining sections of $S,T_0,\ldots,T_{n-1}$, one section of each
degree in $\mathcal G$, and one additional section of degree $F$.
Its cardinality is $n+2+|\mathcal G|$.
For $n=3,4$, the set $\mathcal G$ is explicitly as follows:
\[
\begin{array}{c|l}
 n&\mathcal G\\ \hline
 3&2H-T_1-T_2,\quad 3H-T_1-2T_2,\quad 3H-2T_1-T_2\\[2pt]
 4&2H-T_1-2T_2-T_3,\\
  &3H-T_1-2T_2-2T_3,\quad 3H-2T_1-2T_2-T_3,\\
  &4H-T_1-2T_2-3T_3,\quad 4H-3T_1-2T_2-T_3.
\end{array}
\]
Thus the numbers of minimal generators are eight for $I_3$ and eleven
for $I_4$.
\end{corollary}
\begin{proof}
Since $\chi\geq3$, the numerical bound gives
\[
 \delta(\pi)\leq\frac n4\leq\frac{\chi+2}{2}<\chi.
\]
The cone theorem and Lemma~\ref{lem:MW} therefore apply, and
$\MW(\pi)=0$. The existence and minimality assertions follow from
Theorems~\ref{thm:criterion} and~\ref{thm:bound}. There is precisely one
additional section of degree $F$, because the product of the component
sections spans only a one-dimensional subspace of $H^0(X,\OO_X(F))$.

We verify the two lists directly from the definition of $P$.
For positive integers $p,q$, put
\[
 b_i(p,q)=\min\{ip,(n-i)q\},\qquad
 D_{p,q}=(p+q)H-\sum_{i=1}^{n-1}b_i(p,q)T_i.
\]
For $(n-1)^{-1}p\leq q\leq(n-1)p$, these divisors are nef,
are orthogonal to $S$, and belong to $M^*$: the sequence $b_i$ is
concave, with zero endpoints, and $b_1+b_{n-1}=p+q$.

Conversely, write a nonzero element of $M^*$ as
\[
 A=mH-\sum_{i=1}^{n-1}d_iT_i,\qquad m\geq2,\qquad d_0=d_n=0.
\]
Nefness gives
\[
 2d_i\geq d_{i-1}+d_{i+1},\qquad d_1+d_{n-1}\leq m,
 \qquad
 0\leq d_i\leq\min\{i d_1,(n-i)d_{n-1}\}.
\]
If the $d_i$ are not all zero, their endpoint pair
$(d_1,d_{n-1})$ is positive and lies between the rays $(n-1,1)$ and
$(1,n-1)$. Increasing its smaller coordinate repeatedly preserves this
property and produces a pair $(p,q)$ with $p+q=m$ and
$b_i(p,q)\geq d_i$. If all $d_i$ are zero, take
$(p,q)=(\lfloor m/2\rfloor,\lceil m/2\rceil)$.
In either case, $A-D_{p,q}$ is a nonnegative integral combination of
the nonidentity components.

For $n=3$, the functions $b_i$ are linear throughout the cone between
$(2,1)$ and $(1,2)$; its integral points are generated by
\[
 (2,1),\ (1,1),\ (1,2).
\]
For $n=4$, they are linear on each of the cones between $(3,1)$ and
$(1,1)$, and between $(1,1)$ and $(1,3)$. The integral points in these
cones are generated, respectively, by
\[
 (3,1),\ (2,1),\ (1,1),
 \qquad\text{and}\qquad
 (1,1),\ (1,2),\ (1,3).
\]
These generation statements follow directly because each pair of
consecutive displayed vectors has determinant one and hence is an
integral basis of the lattice in its cone.
Linearity of the $b_i$ then shows that $P$ is generated by the
corresponding $D_{p,q}$ and the nonidentity components.

It remains to check which horizontal generators are indecomposable in
$P$. If one of the displayed degrees were a sum of these generators,
comparison of its fiber degree and the coefficients of $T_1,T_{n-1}$
would give
\[
 \sum_j(p_j+q_j)=p+q,\qquad
 \sum_jp_j\geq p,\qquad \sum_jq_j\geq q.
\]
Thus both inequalities are equalities. Every displayed pair has
$p=1$ or $q=1$, whereas every horizontal summand has both coordinates
positive. Consequently there is exactly one horizontal summand, with
the same pair $(p,q)$; no component summands remain.
These degrees are therefore indecomposable, and they constitute all
of $\mathcal G$. Substitution gives the asserted lists and counts.
\end{proof}

\begin{remark}
The corollary determines the degrees of the generators and proves the
existence of suitable sections. Explicit expressions for these
sections, and equations for the relations among them, require further
work beyond the general criterion.
\end{remark}

\end{document}